\documentclass[12pt,a4paper]{article}
\usepackage[english]{babel}
\usepackage[utf8]{inputenc}
\usepackage{amsmath}
\usepackage{amsfonts}
\usepackage{amssymb}
\usepackage{graphicx}
\usepackage{listings}
\usepackage{listingsutf8}
\usepackage{tikz}
\usepackage{cite}
\usepackage{algorithm2e}
\usepackage{array}
\usepackage{geometry}
\DeclareMathOperator{\rank}{rk}

\DeclareMathOperator{\trace}{Tr}

\begin{document}

\allowdisplaybreaks

\setcounter{tocdepth}{2}

\nocite{*}
\newcounter{currentcounter}

\newcounter{corollarycounter}

\newenvironment{general}[1]{\refstepcounter{currentcounter} \bigskip \noindent \textbf{#1 \thecurrentcounter :}\itshape}{\normalfont \medskip}

\newenvironment{generalbis}[2]{\refstepcounter{currentcounter} \bigskip \noindent \textbf{#1 \thesubsection.\thecurrentcounter ~#2 :}\itshape}{\normalfont \medskip}

\newenvironment{proof}{\noindent \textbf{Proof :} \newline \noindent \hspace*{0.2cm}}{\hspace*{\fill}$\square$ \bigskip}

\newenvironment{proofbis}[1]{\noindent \textbf{#1:} \newline \noindent \hspace*{0.2cm}}{\hspace*{\fill}$\square$ \bigskip}

\newenvironment{famous}[1]{\medskip \noindent \textbf{#1:} \newline \noindent  \itshape}{\normalfont \medskip}

\newenvironment{mycor}{\refstepcounter{corollarycounter} \medskip \noindent \textbf{Corollary \thecorollarycounter :}  \itshape}{\normalfont \medskip}

\newenvironment{court}[1]{\refstepcounter{currentcounter} \smallskip \noindent \textbf{#1 \thesubsection.\thecurrentcounter :} \newline \noindent \itshape}{\normalfont \smallskip}

\newcolumntype{C}{>{$}c<{$}}
\newcommand\toline{\smallskip \newline}
\newcommand\refp[1]{(\ref{#1})}
\renewcommand\mod{\mathrm{~mod~}}
\newcommand\pregcd{\mathrm{gcd}}
\renewcommand\gcd[2]{\pregcd(#1,#2)}
\newcommand\prelcm{\mathrm{lcm}}
\newcommand\lcm[2]{\prelcm(#1,#2)}
\newcommand\naturaliso{\cong}
\newcommand\iso{\simeq}
\newcommand\cross[1]{#1^{\times}}
\newcommand\crosslong[1]{\cross{(#1)}}
\newcommand\poldegree[1]{\mathrm{deg}(#1)}
\newcommand\cardinal[1]{\# \left(#1\right)}
\newcommand\cardinalshort[1]{\##1}
\renewcommand\det{\mathrm{det}}
\newcommand\indicator[1]{[#1]}

\newcommand\sign{\mathrm{sign}}
\newcommand\signature{\mathrm{sgn}}

\newcommand\bb[1]{\mathbb{#1}}
\newcommand\kk{\bb{K}}
\newcommand\ok{\mathcal{O}_{\kk}}
\newcommand\zz{\bb{Z}}
\newcommand\qq{\bb{Q}}
\newcommand\rr{\bb{R}}
\newcommand\cc{\bb{C}}
\newcommand\ff{\bb{F}}
\newcommand\of{\mathcal{O}_{\ff}}
\newcommand\goth[1]{\mathfrak{#1}}
\newcommand\zsz[1]{\zz/#1\zz}
\newcommand\zszcross[1]{\zsz{#1}^{\times}}
\newcommand\dualsimple[1]{#1^{\vee}}
\newcommand\dual[1]{(#1)^{\vee}}

\newcommand\dx[1]{\mathrm{d}#1}
\newcommand\ddx[1]{\frac{\mathrm{d}}{\dx{#1}}}
\newcommand\dkdxk[2]{\frac{\mathrm{d}^{#2}}{\dx{#1}^{#2}}}
\newcommand\partialx[1]{\partial #1}
\newcommand\partialdx[1]{\frac{\partial}{\partialx{#1}}}
\newcommand\partialkdxk[2]{\frac{\partial^{#2}}{\partialx{#1}^{#2}}}
\newcommand\integral[3]{\int_{#1}{#2}\dx{#3}}
\newcommand\integralsimple[1]{\integrale{G}{#1}{\lambda}}
\newcommand\sprod[2]{\langle #1, #2 \rangle}

\newcommand\limi[1]{\lim_{#1 \to + \infty}}

\newcommand\mn[2]{\mathrm{M}_{#1}(#2)}
\newcommand\mnz[1]{\mn{#1}{\zz}}
\newcommand\sln[2]{\mathrm{SL}_{#1}(#2)}
\newcommand\slnz[1]{\sln{#1}{\zz}}
\newcommand\gln[2]{\mathrm{GL}_{#1}(#2)}
\newcommand\glnz[1]{\gln{#1}{\zz}}
\newcommand\resc{\mathrm{resc}}
\newcommand\com{\mathrm{com}}

\newcommand\bars[1]{\underline{#1}}
\newcommand\taubar{\bars{\tau}}
\newcommand\qbar{\bars{q}}
\newcommand\sigmabar{\bars{\sigma}}
\newcommand\rbar{\bars{r}}
\newcommand\nbar{\bars{n}}
\newcommand\mbar{\bars{m}}
\newcommand\mubar{\bars{\mu}}
\newcommand\xbar{\bars{x}}
\newcommand\Xbar{\bars{X}}
\newcommand\abar{\bars{a}}
\newcommand\bbar{\bars{b}}
\newcommand\alphabar{\bars{\alpha}}
\newcommand\betabar{\bars{\beta}}
\newcommand\omegabar{\bars{\omega}}
\newcommand\taubarsj[1]{\taubar^{-}(#1)}
\newcommand\onebar{\bars{1}}

\newcommand\hh{\bb{H}}
\newenvironment{psmallmatrix}
  {\left(\begin{smallmatrix}}
  {\end{smallmatrix}\right)}
  
\newcommand\declarefunction[5]{#1 := \begin{cases}
\hfill #2 \hfill & \to ~ #3 \\
\hfill #4 \hfill & \to ~ #5 \\
\end{cases}}

\newcommand\sume{\sideset{}{_e}\sum}

\newcommand\tendstowhen[1]{\xrightarrow[#1]{}}

\newcommand\pending{$\bigskip \newline \blacktriangle \blacktriangle \blacktriangle \blacktriangle \blacktriangle \blacktriangle \blacktriangle \blacktriangle \bigskip \newline \hfill$}

\newcommand\pendingref{\textbf{[?]}}

\newcommand\opc[1]{\mathcal{O}^{+,\times}_{#1}}
\newcommand\opck{\opc{\kk}}
\newcommand\opcf{\opc{\goth{f}}}
\newcommand\eps{\varepsilon}
\newcommand\units{\ok^{\times}}
\newcommand\unitsf{\of^{\times}}
\newcommand\norm[1]{\mathcal{N}(#1)}

\newcommand\ta{\tilde{a}}
\newcommand\tabar{\tilde{\abar}}
\newcommand\talpha{\tilde{\alpha}}
\newcommand\talphabar{\tilde{\alphabar}}
\newcommand\ts{\tilde{s}}
\newcommand\ttt{\tilde{t}}
\newcommand\tlambda{\tilde{\lambda}}
\newcommand\tc{\tilde{c}}
\newcommand\tcone{\tc_1}
\newcommand\tctwo{\tc_2}
\newcommand\te{\tilde{e}}
\newcommand\teone{\te_1}
\newcommand\tetwo{\te_2}
\newcommand\tmu{\tilde{\mu}}
\newcommand\tA{\tilde{A}}

\newcommand\spvgamma{\Gamma_{\goth{f}, \goth{b}, \goth{a}}(\eps; h)}
\newcommand\spvgammasign{\Gamma_{\goth{f}, \goth{b}, \goth{a}}^{\pm}(\eps; h)}
\newcommand\spvgr{G_{r, \goth{f}, \goth{b}, \goth{a}}^{\pm}(u_1, \dots, u_r; h)}
\newcommand\spvgrrho{G_{r, \goth{f}, \goth{b}, \goth{a}}^{\pm}([\eps_{\rho(1)}|\dots|\eps_{\rho(r)}]; h_\rho)}
\newcommand\spvgrcomplete{I_{r,\goth{f}, \goth{b}, \goth{a}}(\eps_1, \dots, \eps_r ; \bars{h}, \bars{\mu}, \bars{\nu})}

\newcommand\spvgammac{\Gamma_{\goth{f}, \goth{b}, \goth{a}}(\eps; h, \sigma_{\cc})}
\newcommand\spvgrc{G_{r, \goth{f}, \goth{b}, \goth{a}}^{\pm}([\eps_{\rho(1)}|\dots|\eps_{\rho(r)}]; h_\rho, \sigma_{\cc})}
\newcommand\spvgrrhoc{G_{r, \goth{f}, \goth{b}, \goth{a}}^{\pm}([\eps_{\rho(1)}|\dots|\eps_{\rho(r)}]; h_\rho)}
\newcommand\spvgrcompletec{I_{r,\goth{f}, \goth{b}, \goth{a}}(\eps_1, \dots, \eps_r ; \bars{h}, \bars{\mu}, \bars{\nu}, \sigma_{\cc})}

\newcommand\ula{u_{L, \goth{a}}}
\newcommand\ulah{u_{L, \goth{a}, \bars{h}}}
    
\newcommand\hombase{\mathrm{Hom}}
\newcommand\myhom[3]{\hombase_{#1}(#2, #3)}    
\newcommand\homlong[3]{\hombase_{#1}\left(#2, #3\right)}   
\newcommand\homlz{\myhom{\zz}{L}{\zz}}
\newcommand\homlprimez{\myhom{\zz}{L'}{\zz}}
\newcommand\homlc{\myhom{\zz}{L}{\cc}}
\newcommand\homlambdaz{\myhom{\zz}{\Lambda}{\zz}}
\newcommand\homlambdac{\myhom{\zz}{\Lambda}{\cc}}
\newcommand\homvq{\myhom{\qq}{V}{\qq}}
\newcommand\zexc{z}
  
\newcommand\zfone{\mathcal{Z}_{\goth{f}}^1}
  
\newcommand\fracpart[1]{\left\{#1\right\}}
\newcommand\entirepart[1]{\left\lfloor#1\right\rfloor}  
\newcommand\shortsetseparator[1][]{#1|}
\newcommand\setseparator[1][]{~\shortsetseparator[#1]~}
\newcommand\coeff{\mathrm{coeff}}
\newcommand\normalbone[1]{\left(\left(#1\right)\right)}
\newcommand\tensor{\otimes}
\newcommand{\at}[2][]{#1|_{#2}}
\newcommand\omitvar[1]{\widehat{#1}}
\newcommand\kronecker{\delta}
\newcommand\plgt{x}
\newcommand\cohomd{\partial}
\newcommand\cohomdx{\partial^{\times}}
\newcommand\compset[1]{#1^c}
\newcommand\signdet{\sign\,\det}

\newcommand\copen{c^{\circ}}
\newcommand\cclosed{c}
\newcommand\cdual{c^{\vee}}
\newcommand\copendual{c^{\vee, \circ}}
\newcommand\dirac{\delta}

\newcommand\myspan{\mathrm{Span}}
\newcommand\prespanconvex[1]{\mathcal{C}(#1)}
\newcommand\spanconvex{\prespanconvex{V}}
\newcommand\prespancones[1]{\mathcal{K}(#1)}
\newcommand\spancones{\prespancones{V}}
\newcommand\spanconesrr{\prespancones{V_{\rr}}}
\newcommand\spanconesk[1]{\mathcal{K}^{#1}(V)}
\newcommand\spanqcones{\mathcal{K}_{\qq}(V)}
\newcommand\spanqconesrr{\mathcal{K}_{\qq}(V_{\rr})}
\newcommand\prespanwedges[1]{\mathcal{L}(#1)}
\newcommand\spanwedges{\prespanwedges{V}}
\newcommand\spanqwedges{\mathcal{L}_{\qq}(V)}
\newcommand\spanqwedgesrr{\mathcal{L}_{\qq}(V_{\rr})}

\newcommand\coefficient[3]{\mathrm{coeff}[#2^{#3}]\left(#1\right)}
\newcommand\badposition{\mathrm{(BP)}}
\newcommand\symbolparagraph{\S}

\newcommand\alphajk[2]{\alpha^{(#1)}_{#2}}
\newcommand\djk[2]{d^{(#1)}_{#2}}
\newcommand\vjkl{v^{(j)}_{k,l}}

\newcommand\disjointunion{\sqcup}
\newcommand\dkxs{d_k(x, s)}

\begin{center}
\end{center}
\begin{center}
\large \noindent Geometric families of multiple elliptic Gamma functions \linebreak and arithmetic applications, I
\end{center}
\medskip
\begin{center}
Pierre L. L. Morain\footnotemark[1]\footnotetext[1]{Sorbonne Université and Université Paris Cité, CNRS, IMJ-PRG, F-75005 Paris, France. This PhD work is funded by the École polytechnique,  Palaiseau, France.}
\end{center}
\medskip
\begin{center}
\noindent \textbf{Abstract:}
\end{center}
This is the first paper in a series where we study arithmetic applications of the multiple elliptic Gamma functions originated from mathematical physics. The main purpose of this paper is the introduction of a framework for applications of these functions to Hilbert's 12th problem for general number fields with exactly one complex place following recent work by Bergeron, Charollois and Garc\'ia. Namely, we define geometric families of the multiple elliptic Gamma functions, upgrading the construction carried out by Felder, Henriques, Rossi and Zhu for rank $3$ lattices to lattices of higher ranks. These functions enjoy transformation properties under an action of the special linear group $\mathrm{SL}_n(\mathbb{Z})$ for $n \geq 2$ involving some Bernoulli rational functions as their so-called modularity defect. A second purpose of this paper is to use this collection of Bernoulli rational functions to construct $(n-1)$-cocycles for specific subgroups of $\mathrm{SL}_n(\mathbb{Z})$ associated to  groups of units in totally real number fields and use these cocycles to compute partial zeta values at $s=0$.

\normalsize \bigskip \bigskip

\noindent \textbf{Acknowledgments:} 
This work is part of an on-going PhD work and the author would like to thank his advisors Pierre Charollois and Antonin Guilloux at Sorbonne Université for their guidance and their helpful comments on this work.
\newpage

\tableofcontents

\section{Introduction}

This paper deals with the study of multivariate meromorphic functions enjoying transformation properties under an action of the special linear group $\slnz{n}$ on rank $n$ lattices similar to those of $\theta$ functions for $\slnz{2}$. We first consider the following $\theta$ function:
$$ \theta(z, \tau) = \prod_{m \geq 1}(1-e^{2i\pi (m+1)\tau} e^{-2i\pi z})(1-e^{2i\pi m \tau} e^{2i\pi z})$$
It is a holomorphic function of $z \in \cc$ and $\tau \in \hh$ where $\hh = \{t \in \cc \setseparator \Im(t) > 0\}$ is the upper half-plane in $\cc$, satisfying the following relations:
\begin{align*}
\theta(z+1 , \tau) & = \theta(z, \tau) = \theta(z, \tau +1) \\
\theta(z+\tau, \tau) & = -e^{-2i\pi z} \theta(z, \tau) \\
\theta\left(\frac{z}{\tau}, \frac{-1}{\tau}\right) & = \theta(z, \tau) e^{\frac{2i\pi}{\tau} P_2(z, \tau)}
\end{align*}
where the polynomial $P_2(z, \tau)$ is given by:
$$ P_2(z, \tau) = \frac{z^2 + z - z\tau}{2} - \frac{\tau}{4} + \frac{\tau^2-1}{12}  $$
More generally, to any matrix $\gamma = \begin{pmatrix} a & b \\ c & d \end{pmatrix} \in \slnz{2}$ we can associate a rational function $P_{2, \gamma}(z, \tau) \in \qq(z, \tau)$ such that:
$$ \theta\left(\frac{z}{c \tau + d}, \frac{a\tau +b}{c\tau +d}\right) = \theta(z, \tau) e^{2i\pi P_{2, \gamma}(z, \tau)} $$
This may be interpreted as a coboundary relation and shows that $\gamma \to P_{2, \gamma}$ is a $1$-cocycle on $\slnz{2}$ with values in rational functions and that it is split by the $\theta$ function. In this series of articles, we study couples of functions that behave similarly to the couple $(\theta, P_{2, \gamma})$ under an action of a special linear group $\slnz{n}$ of higher dimension $n \geq 3$ on rank $n$ lattices and give an insight on how they might be used to compute invariants in number theory.
 
For $n = 3$, it was shown in the 2000s by Felder and Varchenko \cite{FV} that the elliptic Gamma function introduced by Ruijsenaars \cite{Ruijsenaars} enjoyed modular transformation properties for $\slnz{3}$ similar to those of the $\theta$ function for $\slnz{2}$. This multivariate analytic function is defined on $\cc \times \hh \times \hh$ by:
$$\Gamma(z, \tau, \sigma) = \prod_{m,n \,\geq 0}\left(\frac{1-\exp(2i\pi((m+1)\tau +(n+1)\sigma-z))}{1-\exp(2i\pi(m\tau+n\sigma+z))} \right)$$
and it may be associated with a $2$-cocycle on $\slnz{3}$ with values in rational functions (see section \ref{sectiongeneralgr}). In 2001, Nishizawa \cite{Nishizawa} introduced a whole hierarchy of multiple elliptic Gamma functions which encompass both the $\theta$ function and the elliptic $\Gamma$ function. They are multivariate analytic functions defined for all $r \in \bb{N}$ by:
$$G_r(z, \tau_0, \dots, \tau_r) = \prod_{m_0, \dots, m_r \,\geq 0} \left(1 - e^{2i\pi(-z + \sum_{j = 0}^r (m_j+1)\tau_j)}\right)  \left(1- e^{2i\pi (z + \sum_{j = 0}^r m_j\tau_j)}\right)^{(-1)^r}$$
This definition recovers both the $\theta$ and the elliptic $\Gamma$ functions as $\theta = G_0$ and $\Gamma = G_1$. The $G_r$ functions share similar transformation properties under an action of $\slnz{r+2}$ on the abstract lattice generated by $1, \tau_0, \dots, \tau_r$ (see section \ref{sectiongeneralgr}) and it is natural to study the entire collection obtained when varying $r \in \mathbb{N}$. In particular, when $r \geq 1$, the following pseudo-periodicity relation for the $G_r$ function involves a lower degree $G_{r-1}$ function as:
$$G_r(z + \tau_j, \tau_0, \dots, \tau_r) = G_{r-1}(z, \tau_0, \dots, \omitvar{\tau_j}, \dots, \tau_r)\,G_r(z, \tau_0, \dots, \tau_r)$$
where as usual the notation $\omitvar{\tau_j}$ indicates that the variable $\tau_j$ should be omited. In \cite{FDuke} Felder, Henriques, Rossi and Zhu enriched the theory of the elliptic Gamma function by introducing geometric families of elliptic Gamma functions associated to arbitrary rank $3$ lattices, offering a comprehensive perspective of the underlying geometric phenomena. Namely, if $L$ is a rank $3$ lattice and $a, b \in \homlz $ are two primitive linear forms on $L$, they define general functions $\Gamma_{a,b} : \cc \times \homlc \to \cc$ using Ruijsenaars' elliptic $\Gamma$ function as a building block and show that the collection of functions obtained when varying $a, b$ enjoys nice properties under an action of $\slnz{3}$ on $L$. The two most important properties in that regard are the so-called modular and equivariance properties. For any linearly independent $a, b, c \in \homlz$, there is a rational function $P_{a,b,c} \in \qq[w](\plgt)$ such that the equality:
\begin{equation}\label{introductionmodulargamma}
\Gamma_{a, b}(w, \plgt)\Gamma_{b, c}(w, \plgt)\Gamma_{c, a}(w, \plgt) = \exp(2i\pi P_{a, b, c}(w, \plgt))
\end{equation}
holds for $(w, \plgt)$ in a dense open subset of $\cc \times \homlc \simeq \cc \times \cc^3$.
Additionally, both $\Gamma_{a,b}$ and $P_{a, b, c}$ are \textit{equivariant} under an action of $\slnz{3}$, which means that for all $g \in \slnz{3}$:
\begin{align*}
\Gamma_{g \cdot a, g \cdot b}(w, g \cdot \plgt) & = \Gamma_{a, b}(w, \plgt) \\
P_{g\cdot a, g\cdot b, g \cdot c}(w, g \cdot \plgt) & = P_{a, b, c}(w, \plgt)
\end{align*}
These equivariance properties play a key role in arithmetic applications, as the functions $\Gamma_{a,b}$ and $P_{a,b,c}$ will be evaluated on homology classes associated to specific tori in $\slnz{3}$. 

The construction of the functions $\Gamma_{a,b}$ may already be viewed as the generalisation of a well-known construction for the $\theta$ function. For a rank $2$ lattice $L$, and a linear form $a \in \homlz = \Lambda$, we may also define an equivariant function $\theta_a :  \cc \times \homlc \to \cc$ such that for any pair of linearly independent $a, b \in \Lambda$ there is an equivariant rational function $Q_{a,b} \in \qq[w](\plgt)$ which may be expressed in terms of Dedekind sums and satisfying:
\begin{equation}\label{introductionmodulartheta}
\theta_a(w, \plgt)\theta_b(w, \plgt)^{-1} = \exp(2i\pi Q_{a,b}(w, \plgt))
\end{equation}
In the first part of this article, we upgrade the construction carried out in \cite{FDuke} to higher degree $G_r$ functions (see section \ref{sectiongeneralgr}). For a lattice $\Lambda = \homlz$ of rank $n = r+2 \geq 2$ we introduce a collection of geometric functions $G_{n-2, a_1, \dots, a_{n-1}}$ attached to families of $n-1$ linear forms $a_1, \dots, a_{n-1} \in \Lambda$ (see Proposition \ref{propositiondefgeometricgr}) which are built using Nishizawa's $G_{n-2}$ functions. For rank $n = 2, 3$ lattices, we recover the functions $\theta_a = G_{0, a}(0)$ and $\Gamma_{a,b} = G_{1, a, b}(0)$ respectively. In section \ref{sectiongeneralgr} we show that the collection of $G_{n-2, a_1, \dots, a_{n-1}}$ functions obtained when varying the linear forms $a_1, \dots, a_{n-1}$ satisfy modular and equivariance properties similar to those satisfied by the collections of $\theta_a$ and $\Gamma_{a,b}$ functions. Namely, our first main result is a general version of formulae \refp{introductionmodulargamma} and \refp{introductionmodulartheta} for the $G_{n-2, a_1, \dots, a_{n-1}}$ functions involving a collection of higher degree Bernoulli rational functions $B_{n, a_1, \dots, a_n}$ (see Definition \ref{definitiongeometricbernoulli}). 

\begin{general}{Theorem}\label{theoremmodulargeometricgrintroduction}
Let $L$ be an oriented lattice of rank $n \geq 2$. Let $a_1, \dots, a_n \in \homlz $ be a family of $n$ linearly independent primitive linear forms on $L$. Fix $v \in V/L \simeq \qq^n/\zz^n$ where $V = L \otimes_{\zz} \qq$.
\begin{enumerate}
\item $[\,$Modular property$\,]$ For $(w, x)$ in a dense open subset of $\cc \times \homlc \simeq \cc \times \cc^n$:
\begin{equation}\label{modulargeometricgrintroduction}
\prod_{j= 1}^{n} G_{n-2, a_1, \dots, \omitvar{a_j}, \dots, a_n}(v)(w, x)^{(-1)^{j+1}} = \exp(2i\pi B_{n, a_1, \dots, a_n}(v)(w,x))\end{equation}
\item $[\,$Equivariance relations$\,]$ For any $g \in \slnz{n}$, the following equalities hold in the space of functions $ V/L \times \cc \times \homlc \to \cc$:
\begin{align*}
G_{n-2,g \cdot a_1, \dots, g \cdot a_{n-1}}(g \cdot v)(w, g\cdot x) &= G_{n-2, a_1, \dots, a_{n-1}}(v)(w, x) \\
B_{n, g\cdot a_1, \dots, g\cdot a_n}(g \cdot v)(w, g\cdot x) &= B_{n, a_1, \dots, a_n}(v)(w, x).
\end{align*}
\end{enumerate}
\end{general}

Theorem \ref{theoremmodulargeometricgrintroduction} gives a more comprehensive perspective on the numerous properties of the $G_{n-2}$ functions (see \ref{sectiondefinitiongr}). In particular, formula \refp{modulargeometricgrintroduction} may be understood as a collection of coboundary relations for $\slnz{n}$ as follows. For two sets $A$ and $B$ write $\mathcal{F}(A, B)$ for the set of functions on $A$ with values in $B$. Fix a primitive linear form $a \in \homlz$ as a base point and define the two functions:
$$\psi_{n,a}  := \begin{cases}\slnz{n}^{n-2} & \to \mathcal{F}(V/L \times \cc \times \homlc, \cc) \\ (g_1, \dots, g_{n-2}) & \to \left((v, w, \plgt) \to G_{n-2, a, g_1\cdot a, \dots, (g_1\dots g_{n-2}) \cdot a}(v)(w, \plgt)\right)\end{cases}$$
$$\phi_{n,a}  := \begin{cases} \slnz{n}^{n-1} & \to \mathcal{F}(V/L, \qq[w](\plgt)) \\ (g_1, \dots, g_{n-1}) & \to  B_{n, a, g_1\cdot a, (g_1g_2)\cdot a, \dots, (g_1\dots g_{n-1})\cdot a}(v)(w, \plgt)\end{cases}$$
When $a, g_1\cdot a, \dots, (g_1\dots g_{n-1})\cdot a$ are linearly independent, the first point of Theorem \ref{theoremmodulargeometricgrintroduction} implies that these functions satisfy the multiplicative coboundary relation:
$$\cohomdx\psi_{n,a}(g_1, \dots, g_{n-1})(v)(w, \plgt) = \exp(2i\pi \phi_{n,a}(g_1, \dots, g_{n-1})(v)(w, \plgt)) $$
for any $v \in V/L$ and for $(w, x)$ in a dense open subset of $\cc \times \homlc$. Throughout this series we shall derive arithmetic applications of \textit{both} collections of functions $\psi = (\psi_{n,a})_{n,a}$ and $\phi = (\phi_{n,a})_{n, a}$ to the computation of class field invariants in number fields. 

This series of papers is organized as follows. In this first article we introduce the framework for the whole series and focus on the cocycle properties for the collection of Bernoulli rational functions $B_{n, a_1, \dots, a_n}$ with applications to the computation of partial zeta values at $s=0$ in totally real number fields. In the second article in this series, we show that the modular property \refp{modulargeometricgrintroduction} holds for a wider range of configurations of the linear forms $a_1, \dots, a_n$. We also introduce a smoothing operation for both collections of $G_{n-2, a_1, \dots, a_{n-1}}$ and $B_{n, a_1, \dots, a_n}$ functions which simplifies the right-hand side of \refp{modulargeometricgrintroduction}, turning the collection $\psi$ into a collection of partial $(n-2)$-modular symbols for $\slnz{n}$. Finally, in the third article of this series, we explain how to evaluate these partial modular symbols attached to $G_{n-2}$ functions to  
compute conjectural elliptic units above number fields of degree $n$ with exactly one complex place in the spirit of Hilbert's 12th problem, upgrading the construction carried out for $n=3$ in a recent article by Bergeron, Charollois and Garc\'ia \cite{BCG}. We showcase this with an example for $n=4$. Consider the quartic field $\kk = \qq(z)$ where $z$ is the complex root of the polynomial $x^4 -6x^3-x^2-3x+1$ lying in the upper half-plane. Then the following two quotients of $G_2$ functions evaluated at points in $\kk$:
\begin{align*}\label{exunitintro}
& \frac{G_2\left(\frac{-1}{2}, \frac{-5z^3 + 29z^2 + 15z + 95}{182}, \frac{-6z^3 + 39z^2 - 10z - 47}{182}, \frac{2z^3 - 13z^2 + z - 24}{182}\right)^{-13}}{G_2\left(\frac{-13}{2}, \frac{-5z^3 + 29z^2 + 15z + 95}{14}, \frac{-6z^3 + 39z^2 - 10z - 47}{14}, \frac{2z^3 - 13z^2 + z - 24}{14}\right)^{-1}}, \\ & \frac{G_2\left(\frac{1}{2}, \frac{2z^3 - 13z^2 + z - 24}{182}, \frac{5z^3 - 29z^2 - 15z - 95}{182}, \frac{-2z^3 + 13z^2 + 6z + 143}{182}\right)^{13}}{G_2\left(\frac{13}{2}, \frac{2z^3 - 13z^2 + z - 24}{14}, \frac{5z^3 - 29z^2 - 15z - 95}{14}, \frac{-2z^3 + 13z^2 + 6z + 143}{14}\right)}
\end{align*}
may be computed with high precision and their product is found to be arbitrarly close to the root $r  \approx 4.1210208... - i\cdot5.0617720...$ of the monic polynomial $P = x^8 - 7x^7 + 33x^6 + 49x^5 + 17x^4 + 49x^3 + 33x^2 - 7x + 1$ which is a palindromic polynomial defining an abelian extension of $\kk$. This computation corresponds to the evaluation of a $2$-cocycle built from smoothed $\psi_{4, a}$ functions against a $2$-cycle arising from the group of totally positive units in $\kk$ (see \cite{preprint} for more examples).

The second part of the present paper is devoted to the study of the collection of Bernoulli rational functions $B_{n, a_1, \dots, a_n}$ together with the attached collection of $(n-1)$-cocycles $\phi_{n, a}$, as the restriction of $\phi_{n, a}$ to specific tori will be used in arithmetic applications to compute partial zeta functions in totally real number fields of degree $n$ at $s = 0$. Shifting the focus on the $B_{n, a_1, \dots, a_n}$ functions we get as a consequence of Theorem \ref{theoremmodulargeometricgrintroduction} the additive cocycle relation:
$$\sum_{j = 0}^{n} (-1)^j B_{n, a_0, \dots, \omitvar{a_j}, \dots, a_n}(v)(w, \plgt) \in \zz$$
for any family $a_0, \dots, a_n$ of primitive linear forms in general position, i.e. such that for any $0 \leq j \leq n$, $\rank(a_0, \dots, \omitvar{a_j}, \dots, a_n) = n$. We improve this result in section \ref{sectioncones} by showing that the stronger additive cocycle relation:
\begin{equation}\label{cocyclerelationgeombernintroduction}
\sum_{j = 0}^{n} (-1)^j B_{n, a_0, \dots, \omitvar{a_j}, \dots, a_n}(v)(w, \plgt) = 0
\end{equation}
holds for almost all configurations of $a_0, \dots, a_n$ in the rank $n$ lattice $\Lambda$ (see the bad position condition \refp{badpositioncondition} in section \ref{sectioncones}).
To prove \refp{cocyclerelationgeombernintroduction} we show that the rational functions $B_{n, a_1, \dots, a_n}$ are given by a coefficient extraction in the generating series associated to a rational polyhedral cone and that the cocycle relation they satisfy may be obtained as a specialisation of a cocycle relation satisfied by some indicator functions of closed polyhedral cones. For ease of presentation, we will consider cones in $\qq$-vector spaces but most results we prove on cones may be readily transposed to cones in $\ff$-vector spaces where $\ff$ is any ordered field, say $\rr$ for instance. A polyhedral cone in a $\qq$-vector space $V$ is a set
$$C =  \qq_{\geq 0} v_1 + \dots + \qq_{\geq 0} v_p + \qq_{> 0} v'_1 + \dots + \qq_{> 0}v'_q$$
where the $v_i$'s and the $v'_j$'s are non-zero vectors in $V$. Let us denote by $\spancones$ the $\qq$-algebra generated by the indicator functions of such cones and by $\spanwedges$ the subspace of $\spancones$ generated by the indicator functions of those cones containing some line $\qq v$. We prove a cocycle relation for specific cones which we now define. For non-zero linear forms $a_0, \dots, a_m \in \dualsimple{V}$ define:
$$\cdual(a_0, \dots, a_m)(v) := \begin{cases} 1 \text{ if } \forall\, 0 \leq j \leq m, a_j(v) \geq 0 \\ 0 \text{ otherwise} \end{cases} $$
In section \ref{sectioncones} we prove our second main result which is the technical heart of this paper and might be of independent interest.

\begin{general}{Theorem}\label{theoremcocyclecones}
Let $V$ be a $\qq$-vector space of dimension $n \geq 1$ and let $a_0, \dots, a_n$ be $n+1$ non-zero linear forms on $V$ which generate $\dualsimple{V}$. For all $0 \leq j \leq n$, denote $\eps_j = (-1)^j \signdet(a_0, \dots, \omitvar{a_j}, \dots, a_n) \in \{ -1, 0, 1\}$.
\begin{itemize}
\item[$\mathrm{(i)}$] If there are coefficients $\lambda_0 > 0, \dots, \lambda_n > 0$ such that $\sum_{j = 0}^n \lambda_j a_j = 0$ then the signs $\eps_j$ are all equal to a common sign $\eps$ and 
$$ \sum_{j = 0}^n \eps_j\cdual(a_0, \dots, \omitvar{a_j}, \dots, a_n) \equiv \eps \dirac \mod \spanwedges$$
where $\dirac$ is the Dirac function at $0$.
\item[$\mathrm{(ii)}$] If there is a relation $\sum_{j = 0}^n \lambda_j a_j = 0$ with at least one positive and one negative coefficient among $\lambda_0, \dots, \lambda_n$, then:
$$ \sum_{j = 0}^n \eps_j\cdual(a_0, \dots, \omitvar{a_j}, \dots, a_n) \equiv 0 \mod \spanwedges.$$
\end{itemize}
\end{general}

This theorem may be viewed as a dual theorem to [\hspace{1sp}\cite{CDG}, Theorem 1.1] from which it is inspired. It is interesting to note that for $n=2$ this theorem allows us to recover a result in a recent article by Sharifi and Venkatesh \cite{VS} (see section \ref{sectionconesgenerating}). Our main goal however, is to deduce from this theorem the following corollary on the
cocycle relation \refp{cocyclerelationgeombernintroduction} satisfied by the Bernoulli rational functions (see Proposition \ref{propositionbernoullicocycle} for more details) for almost all configurations of the parameters.

\begin{mycor}\label{corollarybernoulli}
Let $a_0, \dots, a_n$ be $n+1$ linear forms on $V$. Suppose that either $\rank(a_0, \dots, a_n) < n$, or $a_0, \dots, a_n$ generate $\dualsimple{V}$ and there is a relation $\sum_{j = 0}^n \lambda_j a_j = 0$ with at least one positive and one negative coefficient among $\lambda_0, \dots, \lambda_n$. Then:
$$ \sum_{j = 0}^n (-1)^j B_{n, a_0, \dots, \omitvar{a_j}, \dots, a_n}(v)(w,x) = 0.$$
\end{mycor}

In applications the linear forms $a_j$ will be taken of the form $g_j \cdot a$ for some base point $a \in \dualsimple{V}$ and the invertible matrix $g_j$ in some subgroup $U$ of $\slnz{n}$. At the end of section \ref{sectioncones} we give examples of specific subgroups $U$ of $\slnz{n}$ for which any family $g_0, \dots, g_n \in U$ is such that the family $g_0 \cdot a, \dots, g_n \cdot a$ satisfies the hypothesis of Corollary \ref{corollarybernoulli} for any base point $a \in \dualsimple{V}$. The collection of functions $(\phi_{n, a})_{a}$ thus reduces to a collection of homogeneous $(n-1)$-cocycles for $U$. In the arithmetic applications we have in mind to the computation of partial zeta values in totally real number fields at $s=0$, the group $U$ will arise from the group of totally positive units of a given totally real number field. Following \cite{Shintani}, \cite{Colmez}, \cite{Diazydiaz} and \cite{CDG} we will prove that the partial zeta values in totally real number fields may be expressed as combinations of values of the geometric Bernoulli functions $B_{n, a_1, \dots, a_n}$ which appear in the study of the collection of $G_{n-2, a_1, \dots, a_{n-1}}$ functions:

\begin{general}{Theorem}\label{theoremshintani}
Let $\ff$ be a totally real number field of degree $n$. Let $\goth{f}$ be an integral ideal in $\of$ and $\goth{b}$ an integral ideal representing a class in the narrow ray class group $Cl^{+}(\goth{f})$. Then there are explicitly computable cones $c_{\rho} = \cdual(a_{1, \rho}, \dots, a_{n, \rho})$ parametrised by $\rho \in \goth{S}_{n-1}$ and signs $\nu_{\rho} \in \{-1, 0, +1\}$ such that:
$$\zeta_{\goth{f}}(\goth{b}, 0) = \frac1n \sum_{k =1}^n \sum_{\rho \in \goth{S}_{n-1}} \nu_{\rho} B_{n, a_{1, \rho}, \dots, a_{n, \rho}}(1_{\ff})(0, -\sigma_k) $$
where $\sigma_1, \dots, \sigma_n$ are the real embeddings of $\ff$ and the 
$a_{i, \rho}$'s are $\qq$-linear forms on the $n$-dimensional $\qq$-vector space $\ff$. 
\end{general}

We remark that Theorem \ref{theoremshintani} expresses partial zeta values at $s=0$ in a totally real number field $\ff$ as traces of algebraic numbers in $\ff$. These values are therefore rational numbers, as was already known from the theorem of Klingen and Siegel, and already obtained by Shintani using his method. In section \ref{sectionshintani} we will give two explicit examples of such computations for real cubic fields. 

This paper is organized as follows. In section \ref{sectiongeneralgr} we define both collections of $G_{n-2, a_1, \dots, a_{n-1}}$ and $B_{n, a_1, \dots, a_n}$ functions and prove Theorem \ref{theoremmodulargeometricgrintroduction}. In section \ref{sectioncones} we study a cocycle relation for indicator functions of cones and prove Theorem \ref{theoremcocyclecones} which is the technical heart of this paper. We then deduce Corollary \ref{corollarybernoulli} via a coefficient extraction in the generating series attached to cones and show that the $\phi_{n,a}$ functions form a collection of $(n-1)$-cocycles when restricted to specific tori in $\slnz{n}$ arising from unit groups in number fields. In section \ref{sectionshintani} we show that the Bernoulli rational functions $B_{n, a_1, \dots, a_n}$ may be used to express the values of partial zeta functions at $s=0$ in totally real number fields and give the proof of Theorem \ref{theoremshintani}.

\section{Geometric families of multiple elliptic Gamma functions}\label{sectiongeneralgr}
In this section, we recall the properties of the multiple elliptic Gamma functions (the $G_r$ functions) and construct geometric families $G_{r, a_1, \dots, a_{r+1}}$ of these functions attached to a family of $r+1$ linear forms $a_1, \dots, a_{r+1}$, upgrading the $G_r$ functions to collections of equivariant functions for $\slnz{r+2}$ by adapting the construction of the $\Gamma_{a,b}$ functions in \cite{FDuke} to higher degrees.

\subsection{The $G_r$ functions}

\subsubsection{The definition of the $G_r$ functions}\label{sectiondefinitiongr}
We review the definition and properties of the $G_r$ functions inspired by the $\theta$ function of $CM$ theory and the elliptic $\Gamma$ function of Ruijsenaars. The $\theta$ function is defined for $z \in \cc, \tau \in \hh$ by:
$$ \theta(z, \tau) = \prod_{m \geq 1}(1-e^{2i\pi (m+1)\tau} e^{-2i\pi z})(1-e^{2i\pi m \tau} e^{2i\pi z})$$
It is an infinite product which is absolutely convergent and holomorphic over $\cc \times \hh$ enjoying transformation properties under the action of $\slnz{2}$ on the upper half-plane. Indeed, for $z \in \cc, \tau \in \hh$:
\begin{align*}
\theta(z+1 , \tau) & = \theta(z, \tau) = \theta(z, \tau +1) \\
\theta(z+\tau, \tau) & = -e^{-2i\pi z} \theta(z, \tau) \\
\theta\left(\frac{z}{\tau}, \frac{-1}{\tau}\right) & = \theta(z, \tau) e^{\frac{2i\pi}{\tau} P_2(z, \tau)}
\end{align*}
where the polynomial $P_2$ is explicitly given by:
$$ P_2(z, \tau) = \frac{z^2 + z - z\tau}{2} - \frac{\tau}{4} + \frac{\tau^2-1}{12}  $$
The elliptic $\Gamma$ function which was introduced by Ruijsenaars \cite{Ruijsenaars} is defined by:
\begin{equation}\label{defgammar}
\Gamma(z,\tau,\sigma) = \prod_{m,n \geq 0}\left(\frac{1-\exp(2i\pi((m+1)\tau +(n+1)\sigma-z))}{1-\exp(2i\pi(m\tau+n\sigma+z))} \right)
\end{equation}
As pointed out by Spiridonov (see \cite{Spiridonov}) this function was already studied by Jackson in 1905 and implicitly studied in theoretical physics under other names since the 70s. The infinite product \refp{defgammar} is absolutely convergent and holomorphic in $\tau, \sigma \in \hh$ and meromorphic in $z \in \cc$ with poles at points in $\zz + \zz_{\leq 0} \tau + \zz_{\leq 0} \sigma$. In [\hspace{1sp}\cite{FV}, Theorems 3.1 and 4.1] Felder and Varchenko proved that the elliptic $\Gamma$ function enjoys properties under an action of $\slnz{3}$ involving the $\theta$ function. In particular, for any $z, \tau, \sigma$ in the domain of convergence of $\Gamma(z, \tau, \sigma)$:
\begin{align*}
\Gamma(z, \tau, \sigma) & = \Gamma(z, \sigma, \tau)\\
\Gamma(z + 1, \tau, \sigma) = \Gamma(z, \tau + 1, \sigma) & = \Gamma(z, \tau, \sigma+1) = \Gamma(z, \tau, \sigma)  \\
\Gamma(z + \tau, \tau, \sigma) & = \theta(z, \sigma)\Gamma(z, \tau, \sigma) \\
\Gamma(z + \tau + \sigma, \tau, \sigma) & = \Gamma(-z, \tau, \sigma)^{-1} 
\end{align*}
Finally, if $\sigma/\tau \not\in \rr$ then:
\begin{equation} \label{modgamma}
\Gamma(z, \tau, \sigma)^{-1}\Gamma\left(\frac{z}{\tau}, \frac{-1}{\tau}, \frac{\sigma}{\tau}\right) \Gamma\left(\frac{z-\tau}{\sigma}, -\frac{\tau}{\sigma}, -\frac{1}{\sigma}\right)^{-1} = \exp(2i\pi P_3(z,\tau, \sigma))
\end{equation}
where $$P_3(z, \tau, \sigma) = \frac{z^3}{6\tau\sigma} - \frac{\tau + \sigma - 1}{4\tau\sigma}z^2 + \frac{\tau^2 + \sigma^2 + 3\tau\sigma - 3\tau - 3\sigma + 1}{12\tau\sigma}z $$
$$+ \frac{1}{24}(\tau + \sigma - 1)\left(\frac{1}{\tau} + \frac{1}{\sigma} - 1\right)$$ 
The $G_r$ functions introduced by Nishizawa \cite{Nishizawa} using $q$-polylogarithms generalise both the $\theta$ and elliptic $\Gamma$ functions and are defined for $z \in \cc$ and parameters $\tau_0, \dots, \tau_r \in \hh$ by:
$$G_r(z, \taubar) = \prod_{m_0, \dots, m_r \,\geq 0} \left(1 - e^{2i\pi(-z + \sum_{j = 0}^r (m_j+1)\tau_j)}\right)  \left(1- e^{2i\pi (z + \sum_{j = 0}^r m_j\tau_j)}\right)^{(-1)^r}$$
where $\taubar = (\tau_0, \dots, \tau_r)$. This infinite product is absolutely convergent and holomorphic in each parameter $\tau_j \in \hh$ and it is either holomorphic in $z \in \cc$ if $r$ is even or meromorphic in $z \in \cc$ with poles at points of $\zz + \sum_{j = 0}^r \zz_{\leq 0} \tau_j$ if $r$ is odd. For $r = 0, 1$ we recover the definition of the $\theta$ and $\Gamma$ functions so that $G_0 = \theta$ and $G_1 = \Gamma$. The range of the parameters $\tau_j$ can be extended from $\hh$ to $\cc - \rr$ by using nices expressions of the $G_r$ functions as the exponentials of infinite sums involving sinuses and cosines (see [\hspace{1sp}\cite{FV}, formula (15)] and [\hspace{1sp}\cite{Nishizawa}, Proposition 3.6]), namely:
\begin{equation}\label{exponentialformula}
G_r(z, \taubar) = \begin{cases} \exp\left(\sum_{j \geq 1} \frac{1}{(2i)^rj} \frac{\sin(\pi j(2z - (\tau_0+\dots +\tau_r)))}{\prod_{k = 0}^{r} \sin(\pi j \tau_k)}\right) & \text{ if } r \text{ is odd} \\
\exp\left(\sum_{j \geq 1} \frac{2}{(2i)^{r+1}j} \frac{\cos(\pi j(2z - (\tau_0+\dots +\tau_r)))}{\prod_{k = 0}^{r} \sin(\pi j \tau_k)}\right) & \text{ if } r \text{ is even}
\end{cases}
\end{equation}
This expression is valid provided that $\taubar \in (\cc - \rr)^{r+1}$ and $|\Im(2z - (\tau_0+\dots +\tau_r)| < \sum_{j=0}^r |\Im(\tau_j)|$ and allows us to extend the range of parameters to $\taubar \in (\cc - \rr)^{r+1}$ by putting:
\begin{equation}\label{inversionGr}
G_r(z, \tau_0, \dots, \tau_{j-1}, -\tau_j, \tau_{j+1}, \dots, \tau_r) = G_r(z + \tau_j, \tau_0, \dots, \tau_{j-1}, \tau_j, \tau_{j+1}, \dots, \tau_r)^{-1}
\end{equation}
The $G_r$ functions satisfy relations similar to that of the $\theta$ and elliptic $\Gamma$ functions, as proven by Nishizawa in \cite{Nishizawa} and later by Narukawa in \cite{Narukawa}. First, the $G_r$ functions are $1$-periodic in each of their arguments and they satisfy:
\begin{align}
G_r(z + \tau_0 + \dots + \tau_r, \taubar) & = G_r(-z, \taubar)^{(-1)^r} \nonumber \\
G_r(-z, -\taubar) & = G_r(z, \taubar)^{-1} \label{inversiontotale}
\end{align}
Furthermore, if $r \geq 1$, the $G_r$ functions are almost periodic in $z$ with periods $\tau_j$ for $0 \leq j \leq r$ with a correction factor involving a lower degree function:
$$G_r(z + \tau_j, \taubar) = G_{r-1}(z, \tau_0, \dots, \tau_{j-1}, \omitvar{\tau_j}, \tau_{j+1}, \dots, \tau_r)\,G_{r}(z, \taubar)$$
where the notation $\omitvar{\tau_j}$ indicates that the variable $\tau_j$ should be omitted. We now recall the modular property for the $G_r$ functions.

\subsubsection{Modular property and Bernoulli polynomials}\label{sectionordinarybernoulli}
The modular property for the $G_r$ functions was later proved by Narukawa \cite{Narukawa}. To state their theorem, we first need to introduce the multiple Bernoulli polynomials which were implicitly used in the study of the Barnes' multiple $\Gamma$ function \cite{Barnes}. We adopt the following conventions regarding Bernoulli numbers:
\begin{equation}\label{conventionbernoulli}
\frac{t}{e^t-1} = \sum_{k  \geq 0} B_k \frac{t^k}{k!}
\end{equation}
Consider an integer $n \geq 1$. Let $\omegabar = (\omega_1, \dots, \omega_n) \in (\cc-\{0\})^{n}$. We define the multiple Bernoulli polynomials $B_{n,m}^{*}(z, \omegabar)$ with the following generating function:
$$e^{zt}\prod_{j = 1}^n\frac{\omega_j t}{e^{\omega_jt}-1} = \sum_{m \geq 0} B_{n,m}^{*}(z, \omegabar)\frac{t^m}{m!} $$
These polynomials may be expressed explicitly using Bernoulli numbers as:
$$ \frac{1}{m!} B_{n, m}^{*}(z, \omegabar) = \sum_{l = 0}^m \frac{z^l}{l!}\left(\sum_{\substack{k_1 + \dots + k_n = m-l \\ k_j \geq 0}} \left(\prod_{1 \leq j \leq n} \frac{B_{k_j} \omega_j^{k_j}}{k_j!}\right)\right)$$
For $\omegabar \in (\cc-\{0\})^{n},~ B_{n,m}^{*}(z, \omegabar)$ is a degree $m$ homogeneous polynomial in $n+1$ variables with rational coefficients, which is symmetric in the $n$ variables of $\omegabar$. In \cite{Narukawa}, Narukawa used the rescaled homogeneous rational functions $B_{n,m}(z, \omegabar) = (\prod_{j = 1}^{n}\omega_j^{-1})B_{n, m}^{*}(z, \omegabar)$. These obey many relations which can easily be obtained from the properties of the generating function (see \cite{Narukawa}).
We will be most interested by the diagonal polynomials $B_{n, n}^*$. For instance, the polynomial $B_{2,2}^*$ is given by:
$$(\omega_1\omega_2)B_{2,2}(z, \omega_1, \omega_2) = B_{2,2}^{*}(z, \omega_1, \omega_2) = z^2 - z(\omega_1 + \omega_2) + \frac{\omega_1^2+\omega_2^2+3\omega_1\omega_2}{6}.$$
Narukawa's theorem  (see [\hspace{1sp}\cite{Narukawa}, Theorem 7]) can be stated as follows. Let $n \geq 2$. Fix $\omegabar \in (\cc - \{0\})^{n}$ and suppose that $\omega_j/\omega_k \in \cc - \rr$ for all $1 \leq j \neq k \leq n$. Then for $z \in \cc$ outside the discrete set of poles of the left-hand side:
\begin{equation} \label{modGr}
\prod_{j = 1}^{n} G_{n-2}\left(\frac{z}{\omega_j}, \left(\frac{\omega_k}{\omega_j}\right)_{k \neq j}\right) = \exp\left(\frac{-2i\pi}{n!}B_{n,n}(z, \omegabar)\right)~.
\end{equation} 

\subsubsection{Distribution relations}

The arithmetic applications we have in mind for the $G_r$ functions are inspired by the construction of Siegel units using the $\theta = G_0$ function and the construction of associated cocycles (see for instance \cite{DarmonPozziVonk}). To extend the analogy, we end this section by proving that the $G_r$ functions satisfy distribution relations similar to the distribution relations satisfied by Siegel units.

\begin{general}{Proposition}\label{distributiongr}
Consider an integer $N \geq 2$. Then the following distribution relations hold:
\begin{align*}
\prod_{k = 0}^{N-1} G_r\left(z + \frac{k}{N}, \taubar\right) & = G_r(Nz, N\taubar)\\ \prod_{k = 0}^{N-1} G_r\left(z + \frac{k}{N}\tau_l, \tau_0, \dots, \tau_l, \dots, \tau_r\right) & = G_r\left(z, \tau_0, \dots, \frac{\tau_l}{N}, \dots, \tau_r\right) 
\end{align*}
These two relations together give the complete distribution relation:
\begin{equation*}
\prod_{k, k_0, \dots, k_r = 0}^{N-1} G_r\left(z + \frac{k + k_0\tau_0 + \dots + k_r\tau_r}{N}, \tau_0, \dots, \tau_r\right) =  G_r(Nz, \tau_0, \dots, \tau_r)
\end{equation*}
\end{general}

\begin{proof}
The first relation follows from the standard cyclotomic relation:
$$\prod_{k = 0}^{N-1} (1 - e^{2i\pi k/N} y) = 1 - y^N $$
applied here to both $y_0 = -z + \sum_{j = 0}^r (m_j+1)\tau_j$ and $y_1 = z + \sum_{j = 0}^r m_j\tau_j $ for all $m_0, \dots, m_r \geq 0$.

The second relation is obtained by straightforward computation using the set identities:
$$\zz_{>0} - \{0, 1/N, \dots, (N-1)/N\} = \zz_{>0}/N,~~\zz_{\geq 0} + \{0, 1/N, \dots, (N-1)/N\} = \zz_{\geq 0}/N.$$
\end{proof}

\subsection{Geometric $G_{r, \abar}$ functions}\label{sectiongeometricgr}

In this section we upgrade the construction carried out by Felder, Henriques, Rossi and Zhu to higher degree and define geometric families $G_{r, a_1, \dots, a_{r+1}}$ parametrised by $r+1$ linear forms of the $G_r$ functions that encompass the geometric families $\Gamma_{a,b}$ of the elliptic Gamma function. We then prove Theorem \ref{theoremmodulargeometricgrintroduction}, i.e. the modular property for the $G_{r, a_1, \dots, a_{r+1}}$ functions and their equivariance property under the action of $\slnz{r+2}$. The construction of our $G_{r, a_1, \dots, a_{r+1}}$ functions is adapted from the construction of the elliptic $\Gamma_{a,b}$ functions as we use a generalised version of the alternative definition given by [\hspace{1sp}\cite{FDuke}, Proposition 3.5] and then revert the computations to prove that the definition is indeed valid. 

We start by giving a precise geometric setup. Consider a $\qq$-vector space $V$ of dimension $n$ and a rank $n$ lattice $L \subset V$. Denote $\Lambda := \homlz \subset \dualsimple{V} = \myhom{\qq}{V}{\qq}$. The lattice $L$ is then canonically isomorphic to $\homlambdaz$. Fix $B = [e_1, \dots, e_n]$ a $\zz$-basis of $L$ and denote $C = [f_1, \dots, f_n]$ the dual $\zz$-basis of $\Lambda$ such that for all $1 \leq j, k \leq n$, $f_j(e_k) = \kronecker_{jk}$ where $\kronecker_{jk}$ is Kronecker's symbol. This fixes orientation forms on $L$ and $\Lambda$ given by $\det_B$ and $\det_C$ respectively. When there is no risk of confusion we drop the subscripts and write as usual $\det$ for the orientation forms. This also fixes an action of $\slnz{n}$ on $L$ by left multiplication, i.e. $g \cdot \alpha = g\alpha$ and the contragredient action of $\slnz{n}$ on $\Lambda$ given by $g \cdot a = ag^{-1}$. In particular, the pairing:
$$ \begin{cases} \Lambda \times L &\to \zz \\ (a, \alpha) & \to a(\alpha) \end{cases}$$
is equivariant under the action of $\slnz{n}$, which means that for any $g \in \slnz{n}$ and any $(a, \alpha) \in \Lambda \times L$, $(g\cdot a)(g \cdot \alpha) = a(\alpha)$. Note that the basis $C$ is also a basis of the $\cc$-vector space $\homlc \simeq \cc^n$ which induces an action of $\slnz{n}$ on $\homlc$ extending the action of $\slnz{n}$ on $\Lambda$. We will now consider families of linearly independent primitive linear forms in $\Lambda$ and for $r = n-2$ we define functions $G_{r, a_1, \dots, a_{r+1}}$ attached to theses families. Recall that an element $w$ in a $\zz$-module $\Lambda$ is called primitive if for all $(n, w') \in \zz_{\geq 1} \times \Lambda$, $w = nw' \Rightarrow n = 1$. In order to define the functions $G_{r, a_1, \dots, a_{r+1}}$ properly, we now define positive dual families which will be used in the rest of this work.

\begin{general}{Definition}\label{definitionpositivedualfamily}
Let $(a_1, \dots, a_m)$ be a family of $m$ linearly independent elements in a lattice $\Lambda \simeq \homlz$. We call $(\alpha_1, \dots, \alpha_m) \in L$ a positive dual family to $\abar = (a_1, \dots, a_m)$ if for all $1 \leq j \leq m$ the following holds:
$$ a_j(\alpha_j) > 0, ~~~~a_k(\alpha_j) = 0, \, \forall\, k \neq j$$
\end{general}

\noindent The two important cases in this work will be those where $m = n-1$ and $m = n$ in the lattice $\Lambda$ of rank $n$. The following lemma shows that in these cases two positive dual families to the same family $\abar$ are closely related.

\begin{general}{Lemma}\label{lemmapositivedualfamily}
Let $\Lambda$ be a lattice of rank $n$ with an orientation form $\det$.
\begin{itemize}
\item[$\mathrm{(i)}$] Let $\abar = (a_1, \dots, a_{n-1})$ be a family of $n-1$ linearly independent elements in $\Lambda$. If $(\alpha_1, \dots, \alpha_{n-1})$ and $(\alpha'_1, \dots, \alpha'_{n-1})$ are two positive dual families to $\abar$ in $L = \homlambdaz$ then there are rational numbers $t_1, \dots, t_{n-1}$ such that:
$$a_j(\alpha'_j) \alpha_j = a_j(\alpha_j)\alpha'_j + t_j \det(a_1, \dots, a_{n-1}, \cdot)$$
\item[$\mathrm{(ii)}$] A family $\abar = (a_1, \dots, a_n)$ of $n$ linearly independent elements in $\Lambda$ has exactly one positive dual family $\alphabar \in L = \homlambdaz$ containing only primitive vectors. We call this family the primitive positive dual family to $\abar$. Any other positive dual family $\beta_1, \dots, \beta_n$ to $\abar$ satisfies for all $1 \leq j \leq n, \beta_j = m_j \alpha_j$ for some integer $m_j > 0$.
\end{itemize}
\end{general}

\begin{proof}
(i) For any $1 \leq j \leq n-1$, consider $\gamma_j = a_j(\alpha'_j)\alpha_j-a_j(\alpha_j)\alpha'_j$. Then for all $1 \leq k \leq n-1$, we get $\gamma_j(a_k) = 0$. This means that either $\gamma_j = 0$ and then $t_j = 0$ or the $\qq$-linear forms $\gamma_j$ and $\det(a_1, \dots, a_{n-1}, \cdot)$ defined on the $\qq$-vector space $\dualsimple{V}$ share the same kernel, and therefore, they must be linearly dependent. This proves the first claim.

(ii) Suppose that $\alphabar, \,\alphabar'$ are two positive dual families to $a_1, \dots, a_n$ in $L$ consisting of primitive vectors. Consider once again $\gamma_j =  a_j(\alpha'_j)\alpha_j-a_j(\alpha_j)\alpha'_j$ such that for all $1 \leq k \leq n$, $\gamma_j(a_k) = 0$. The family $(a_1, \dots, a_n)$ is a basis of the $\qq$-vector space $\dualsimple{V}$, so $\gamma_j = 0$ for all $1 \leq j \leq n$. This means that $a_j(\alpha'_j)\alpha_j = a_j(\alpha_j)\alpha'_j$ for all $1 \leq j \leq n$. There are only two opposite primitive vectors in the line $\qq \alpha_j$ so if both $\alpha_j$ and $\alpha'_j$ are primitive it must be that $\alpha_j = \alpha'_j$ because both lie in the half-plane $\{v \in V \setseparator a_j(v) > 0\}$. Thus, if such a family exists, it is unique. Regarding the existence, we may write $a_j = \sum_{k = 1}^{n} a_{j,k} f_k$ and define the matrix:
$$ A = \begin{pmatrix} a_{1,1} & \dots & a_{1, n} \\ \vdots & \vdots & \vdots \\ a_{n, 1} & \dots & a_{n, n}\end{pmatrix} $$
If $\epsilon$ is the sign of $\det(A)$ then the matrix $B = \epsilon.\mathrm{com}(A)^{T}$ defines a positive dual family to $\abar$ which we use to compute the primitive positive dual family to $\abar$. Define $\Delta_{i,j} = (-1)^{(i+j)}\det(a_{u,v})_{u \neq i, v \neq j}$ and $\delta_j = \pregcd(\Delta_{1, j}, \dots, \Delta_{n, j})$. Then the primitive positive dual family $\alphabar$ to $\abar$ is given explicitly by: 
$$\alpha_k = \sum_{j = 1}^n \frac{\epsilon\Delta_{k,j}}{\delta_j} e_k $$
If $\beta_1, \dots, \beta_n$ is any other positive dual family to $\abar$ then there are unique integers $m_1, \dots, m_n > 0$ such that $\beta_1/m_1, \dots, \beta_n/m_n$ are primitive vectors in $L$. The family $\beta_1/m_1, \dots, \beta_n/m_n$ is a primitive positive dual family to $\abar$ in $L$, therefore it must be equal to $\alpha_1, \dots, \alpha_n$. This shows that $\beta_j = \alpha_j m_j$ for all $1 \leq j \leq n$.  
\end{proof}

We may now define the geometric variants of the $G_r$ functions using lemma \ref{lemmapositivedualfamily} and adapting [\hspace{1sp}\cite{FDuke}, Proposition 3.5] to higher rank lattices:

\begin{general}{Proposition}\label{propositiondefgeometricgr}
Let $n \geq 2$ and set $r = n-2$. Let $\abar = (a_1, \dots, a_{r+1})$ be a family of $r+1$ linearly independent primitive vectors in the oriented lattice $\Lambda$ of rank $n = r+2$. There is a unique integer $s > 0$ and a unique primitive element $\gamma \in L = \homlambdaz$ such that $\det(a_1, \dots, a_{r+1}, \cdot) = s\gamma$. Fix a vector $v \in V/L$. For any choice of positive dual family $\alphabar = (\alpha_1, \dots, \alpha_{r+1})$ to $\abar$ in $L$ the function defined by the finite product:
\begin{equation}\label{defgeom}
G_{r, \abar}^{\alphabar}(v)(w, \plgt) := \prod_{\delta \in F(\abar, \alphabar, v)/\zz\gamma} G_r\left(\frac{w+\plgt(\delta)}{\plgt(\gamma)}, \frac{1}{\plgt(\gamma)}\plgt(\alphabar) \right)
\end{equation}
where 
$$F(\abar, \alphabar, v) = \{ \delta \in v+ L \setseparator \forall\,1 \leq j \leq r+1,~0 \leq a_j(\delta) < a_j(\alpha_j) \}$$
is well defined for $(w, x)$ in a dense open set of the $\cc$-vector space $\cc \times \homlc \iso \cc \times \cc^{r+2}$. Furthermore, it is independent of the choice of $\alphabar$.
\end{general}

\begin{proof} 
Consider the following open subset of $\homlc$:
$$U(\abar, \alphabar) = \{ \plgt \in \homlc \setseparator \plgt(\gamma) \neq 0,~ \plgt(\alpha_j)/\plgt(\gamma) \not\in \rr, \forall\, 1 \leq j \leq r+1\} $$
The set $U(\abar, \alphabar)$ is the complementary set of a finite union of $\rr$-vector spaces of dimension $\leq 2r+3$ in the $\rr$-vector space $\homlc \simeq \rr^{2r+4}$, therefore it is a dense open subset of $\homlc$ endowed with the finite dimensional $\rr$-vector space topology. For any $\plgt \in U(\abar, \alphabar)$, define:
$$ S(\abar, \alphabar, v, \plgt) = \bigcap_{\delta \in F(\abar, \alphabar, v)} \left\{ w \in \cc \setseparator[\big] \substack{w + \plgt(\delta) \not\in \zz\plgt(\gamma) + \sum_{j = 1}^{r+1} \zz_{> 0} \plgt(\alpha_j)\\ w + \plgt(\delta) \not\in \zz\plgt(\gamma) + \sum_{j = 1}^{r+1} \zz_{\leq 0} \plgt(\alpha_j)}  \right\}$$ 
Then the right-hand side of \refp{defgeom} is well-defined and non-zero on:
\begin{equation}\label{definitionomega}
\Omega(\abar, \alphabar, v) = \bigcup_{\plgt \in U(\abar, \alphabar)} \{(w, \plgt) \setseparator w \in S(\abar, \alphabar, v, \plgt)\}
\end{equation}
which is a dense open subset of $\cc \times \homlc$ endowed with the finite dimensional $\cc$-vector space topology.

Let us now prove that the right-hand side of \refp{defgeom} is indeed independent of the choice of positive dual family $\alphabar$. To achieve this, we will show that for any other choice of positive dual family $\alphabar'$ to $\abar$, the functions $G_{r, \abar}^{\alphabar}$ and $G_{r, \abar}^{\alphabar'}$ coincide. Let us first compute the right-hand side of \refp{defgeom} explicitly using the definition of the $G_r$ function. Consider $(w,\plgt) \in \Omega(\abar, \alphabar, v)$. Put for $1 \leq j \leq r+1$, $d_j = \pm 1$ such that $d_j\plgt(\alpha_j)/\plgt(\gamma) \in \hh$. Put also $D = \sum_{j = 1}^{r+1} (d_j-1)/2$. Then using the inversion relation \refp{inversionGr} we get:
\begin{multline*}
G_{r, \abar}^{\alphabar}(v)(w, \plgt)^{(-1)^D}  = \\
\prod_{\delta \in F(\abar, \alphabar, v)/\zz\gamma} G_r\left(\frac{w+\plgt(\delta)}{\plgt(\gamma)} + \sum_{j = 1}^{r+1} \frac{d_j-1}{2}\frac{\plgt(\alpha_j)}{\plgt(\gamma)}, \left(\frac{1}{\plgt(\gamma)}d_j\plgt(\alpha_j)\right)_{1 \leq j\leq r+1}\right)
\end{multline*}
which by definition of the ordinary $G_r$ functions is given by:
\begin{multline*}
G_{r, \abar}^{\alphabar}(v)(w, \plgt)^{(-1)^D} = \prod_{\delta \in F(\abar, \alphabar, v)/\zz\gamma} \prod_{\mbar\,\geq 0} \Big[\left(1-e^{2i\pi\left(\sum_{j = 1}^{r+1} \frac{(d_jm_j+(1+d_j)/2)\plgt(\alpha_j)}{\plgt(\gamma)} - \frac{w + \plgt(\delta)}{\plgt(\gamma)}\right)}\right) \\
\times \left(1-e^{2i\pi\left(\sum_{j = 1}^{r+1} \frac{(d_jm_j+(d_j-1)/2)\plgt(\alpha_j)}{\plgt(\gamma)} + \frac{w + \plgt(\delta)}{\plgt(\gamma)}\right)}\right)^{(-1)^r}\Big]
\end{multline*}
Let us denote by  $C^{+}(\abar, \alphabar, v, \plgt)$ the set of $\delta' \in v + L$ satisfying for all $1 \leq j \leq r+1$:
$$\begin{cases}
a_j(\delta') \geq 0, \text{ if } d_j = 1 \\
a_j(\delta') < 0, \text{ if } d_j = - 1
\end{cases}$$
and similarly denote by $C^{-}(\abar, \alphabar, v, \plgt)$ the set of $\delta' \in v + L$ satisfying for all $1 \leq j \leq r+1$:
$$\begin{cases}
a_j(\delta') \geq 0, \text{ if } d_j = -1 \\
a_j(\delta') < 0, \text{ if } d_j =  1.
\end{cases}$$
Consider $\delta' \in C^{+}(\abar, \alphabar, v, \plgt)$. If $d_j =1$ then $a_j(\delta') \geq 0$ so that performing Euclidian division by $a_j(\alpha_j)$ gives a unique integer $m_j \geq 0$ satisfying $0 \leq a_j(\delta' - m_j \alpha_j) < a_j(\alpha_j)$. On the contrary, if $d_j = -1$ then $a_j(\delta') < 0$ and there exists a unique integer $m_j > 0$ such that $0 \leq a_j(\delta' + m_j\alpha_j) < a_j(\alpha_j)$. Then for all $1 \leq k \leq r+1$: 
$$0 \leq a_k\left(\delta' - \sum_{j = 1}^{r+1} d_jm_j\alpha_j \right) < a_k(\alpha_k)$$
and thus $\delta' - \sum_{j = 1}^{r+1} d_jm_j\alpha_j \in F(\abar, \alphabar, v)$. This shows that the cone $C^{+}(\abar, \alphabar, v, \plgt)$ can be written as a disjoint union:
$$C^{+}(\abar, \alphabar, v, \plgt) = \bigcup_{\delta \in F(\abar, \alphabar, v)} \bigcup_{\mbar\, \geq 0} \left\{\delta + \sum_{j = 1}^{r+1} (d_jm_j + (d_j-1)/2)\alpha_j\right\}.$$
A similar argument applied to $C^{-}(\abar, \alphabar, v, \plgt)$ gives the decomposition:
$$C^{-}(\abar, \alphabar, v, \plgt) = \bigcup_{\delta \in F(\abar, \alphabar, v)} \bigcup_{\mbar\, \geq 0} \left\{\delta - \sum_{j = 1}^{r+1} (d_jm_j + (d_j+1)/2)\alpha_j\right\}.$$
Thus, the expression $G_{r, \abar}^{\alphabar}(v)(w, \plgt)^{(-1)^D}$ is equal to:
$$ \prod_{\delta' \in C^{-}(\abar, \alphabar, v, \plgt)/\zz\gamma}\left(1-e^{-2i\pi\left(\frac{w+\plgt(\delta')}{\plgt(\gamma)}\right)}\right)\prod_{\delta' \in C^{+}(\abar, \alphabar, v, \plgt)/\zz\gamma}\left(1-e^{2i\pi\left(\frac{w+\plgt(\delta')}{\plgt(\gamma)}\right)}\right)^{(-1)^r}$$
Then, we only need to show that the sets $C^{\pm}(\abar, \alphabar, v, \plgt)$ are independent of the choice for $\alphabar$. Consider another positive dual family $\alphabar'$ to $\abar$. Write $a_j(\alpha_j) = s_j > 0$ and $a_j(\alpha'_j) = s'_j>0$. Then from lemma \ref{lemmapositivedualfamily} there is a rational number $t_j$ such that:
$$s_j\alpha'_j = s'_j\alpha_j + t_j\gamma$$
which gives
$$d_js_j\frac{\plgt(\alpha'_j)}{\plgt(\gamma)} = d_js'_j\frac{\plgt(\alpha_j)}{\plgt(\gamma)} + d_jt_j\in \hh$$
This shows that the signs $d_j$ (and therefore also $D$) are independent of the choice for $\alphabar$. By construction we get $C^{\pm}(\abar, \alphabar, v, \plgt) = C^{\pm}(\abar, \alphabar', v, \plgt)$ and the definition of the geometric $G_{r, \abar}^{\alphabar}$ function is independent of the choice for $\alphabar$. Notice that the sets $U(\abar) = U(\abar, \alphabar), S(\abar, v, \plgt) = S(\abar, \alphabar, v, \plgt), \Omega(\abar, v) = \Omega(\abar, \alphabar, v)$ are all independent of the choice for $\alphabar$.
\end{proof}

From now on we denote by $G_{r, \abar} := G_{r, \abar}^{\alpha}$ the geometric $G_r$ function associated to $\abar$ for any suitable choice of $\alphabar$. By convention, when $a_1, \dots, a_{r+1}$ are not linearly independent, we define $G_{r, \abar}$ to be the constant function equal to $1$. When $r = 0$ and $r = 1$ we recover the geometric variants of the $\theta$ and elliptic $\Gamma$ function as:
$$\theta_a(w, \plgt) = G_{0, a}(0)(w, \plgt), ~~ \Gamma_{a,b}(w, \plgt) = G_{1, a, b}(0)(w, \plgt).$$
Thus, we will often write $\theta_a(v)(w, \plgt)$ for $G_{0, a}(v)(w, \plgt)$ as well as $\Gamma_{a,b}(v)(w,\plgt)$ for $G_{1, a, b}(v)(w, \plgt)$. We remark that Proposition \ref{propositiondefgeometricgr} gives a definition of our $G_{r, a_1, \dots, a_{r+1}}$ functions as $G_r$ functions for specific cones in a rank $r+2$ lattice. This construction may be compared to Winding's construction of $G_r$ functions attached to cones in a rank $r+1$ lattice \cite{Winding}. We argue that we may think of our $G_{r, a_1, \dots, a_{r+1}}$ as objects in projective geometry whereas Winding's construction belongs to affine geometry and both should be related in some sense.

To express the transformation properties of the function $G_{r, \abar}$ under the action of $\slnz{r+2}$ we need to introduce the family of geometric Bernoulli rational functions which encompass both $Q_{a,b}$ and $P_{a,b,c}$ appearing in formulae \refp{introductionmodulartheta} and \refp{introductionmodulargamma} respectively. Recall that we have defined the polynomials $B_{n, n}^{*}$ in section \ref{sectionordinarybernoulli}.

\begin{general}{Definition}\label{definitiongeometricbernoulli}
Let $a_1, \dots, a_n$ be a family of $n$ linearly independent primitive vectors in the oriented lattice $\Lambda$ of rank $n$ and $\alpha_1, \dots, \alpha_n$ be its primitive positive dual family in $L$ (see lemma \ref{lemmapositivedualfamily}). Let $\epsilon$ be the sign of $(-1)^{n}\det(a_1, \dots, a_n)$. We define the geometric Bernoulli polynomial attached to $a_1, \dots, a_n$ and to $v \in V/L$ on $\cc \times \homlc \simeq \cc \times \cc^n$ by the finite sum:
$$ B^{*}_{n, a_1, \dots, a_n}(v)(w, \plgt) := \frac{\epsilon}{n!} \sum_{\delta \in F(\abar, v)} B^{*}_{n,n}(w + \plgt(\delta), \plgt(\alpha_1), \dots, \plgt(\alpha_{n}))$$
where
$$F(\abar, v) = \{\delta \in v + L, 0 \leq a_j(\delta) < a_j(\alpha_j), \forall\,1 \leq j \leq n \} $$ 
is a finite set. The geometric Bernoulli polynomial $B^{*}_{n, a_1, \dots, a_n}(v)$ is a degree $n$ homogeneous polynomial in $n+1$ variables on $\cc \times \homlc \simeq \cc \times \cc^{n}$, with rational coefficients depending on $a_1, \dots, a_n$ and $v$.
\end{general}

Here we give a simple example and recover the Bernoulli polynomials $B_{n, n}^{*}$ from section \ref{sectionordinarybernoulli}. Suppose that the $a_j$ are the vectors of the basis $C$, which means in coordinates that $\forall\, 1 \leq j \leq n,\,a_j = (0, \dots, 0, 1, 0, \dots, 0)$ with a $1$ at the $j$-th position. The corresponding primitive positive dual family is given by the vectors of the basis $B$, which in coordinates gives $\alpha_j =  (0, \dots, 0, 1, 0, \dots, 0)^{T}$ with a $1$ at the $j$-th position. Suppose further $v = 0$, so that $F = F(\abar, 0) = \{0\}$.  Denote by $\tau_j$ the value of $x$ on $\alpha_j$. Then:
$$ B^{*}_{n, a_1, \dots, a_n}(0)(w, \plgt) := \frac{(-1)^n}{n!}B^{*}_{n,n}(w, \tau_1, \dots, \tau_n).$$

\noindent The modularity property for the geometric $G_{n-2, \abar}$ functions will involve the rescaled degree $0$ homogeneous rational functions
\begin{equation}\label{definitionrescaledbernoulli}
B_{n, a_1, \dots, a_n}(v)(w, \plgt) = \left(\prod_{j = 1}^{n} \plgt(\alpha_j)\right)^{-1}B^{*}_{n, a_1, \dots, a_n}(v)(w, \plgt).
\end{equation}
If $a_1, \dots, a_n$ are linearly dependent, we set by convention $B_{n, \abar} = 0$. For $n=2$ and $n=3$ the rational functions $B_{2, a, b}(0)$ and $B_{3, a, b, c}(0)$ may be identified with the rational functions $Q_{a, b}$ and $P_{a,b,c}$ appearing in formulae \refp{introductionmodulartheta} and \refp{introductionmodulargamma} respectively and we may now generalise these formulae by proving Theorem \ref{theoremmodulargeometricgrintroduction}. \bigskip

\begin{proofbis}{Proof of Theorem \ref{theoremmodulargeometricgrintroduction}}
Let $a_1, \dots, a_n$ be a family of $n$ linearly independent elements in the rank $n$ lattice $\Lambda$ and let $\alpha_1, \dots, \alpha_n$ be the primitive positive dual family to $a_1, \dots, a_n$ in $L = \homlambdaz$. Fix $v \in V/L$.

1. We first show that the modular property for the geometric $G_{n-2, \abar}$ functions is a consequence of the modular property \refp{modGr} for the ordinary $G_{n-2}$ function. Write $\epsilon$ for the sign of $(-1)^n\det(a_1, \dots, a_n)$. Then for all $1 \leq j \leq n$ there is a positive integer $s_j$ such that:
$$\det(a_1, \dots, \omitvar{a_j}, \dots, a_n, \cdot) = s_j.\epsilon.(-1)^j\alpha_j$$
Then, by definition, for any $1 \leq j \leq n$:
$$G_{n-2, ((a_k)_{k \neq j})}(v)(w, \plgt)  = \prod_{\delta \in F_j/\zz\alpha_j}G_{n-2}\left(\frac{w + \plgt(\delta)}{\epsilon(-1)^{j}\plgt(\alpha_j)}, \left(\frac{\plgt(\alpha_k)}{\epsilon(-1)^{j}\plgt(\alpha_j)}\right)_{k \neq j}\right)$$
where
$$F_j = \{\delta \in v + L, 0 \leq a_k(\delta) < a_k(\alpha_k), \forall\, 1 \leq k \neq j \leq n \} $$ 
and using the inversion relation \refp{inversiontotale} we get:
$$G_{n-2, ((a_k)_{k \neq j})}(v)(w, \plgt) = \prod_{\delta \in F_j/\zz\alpha_j}G_{n-2}\left(\frac{w + \plgt(\delta)}{\plgt(\alpha_j)}, \left(\frac{\plgt(\alpha_k)}{\plgt(\alpha_j)}\right)_{k \neq j}\right)^{(-1)^{j}.\epsilon}.$$
This equality holds in the dense open subset $\Omega_j = \Omega(a_1, \dots, \omitvar{a_j}, \dots, a_n ; v)$ of $\cc \times \homlc$ (see formula \refp{definitionomega}).
Consequently, the following equality holds on the dense open subset $\Omega = \cap_{j = 1}^n \Omega_j$ of $\cc \times \homlc$ for the whole product:
$$\prod_{j= 1}^{n} G_{r, ((a_k)_{k \neq j})}(v)(w, \plgt)^{(-1)^{j+1}} = \prod_{j= 1}^n  \prod_{\delta \in F_j/\zz\alpha_j}G_{n-2}\left(\frac{w + \plgt(\delta)}{\plgt(\alpha_j)}, \left(\frac{\plgt(\alpha_k)}{\plgt(\alpha_j)}\right)_{k \neq j}\right)^{-\epsilon}.$$
Put $F = \{\delta \in v + L, 0 \leq a_k(\delta) < a_k(\alpha_k), \forall\,1 \leq k \leq n\}$. Then we may write uniformly $F \iso F_j/\zz\alpha_j$ for all $1 \leq j \leq n$. Using Narukawa's theorem (see \refp{modGr}) for each $\delta$ in the finite set $F$ yields:
$$\prod_{j= 1}^{n} G_{n-2, ((a_k)_{k \neq j})}(v)(w, \plgt)^{(-1)^{j+1}} = \prod_{\delta \in F}
\exp\left(\frac{2i\pi\epsilon}{n!}B_{n, n}(w+\plgt(\delta), \plgt(\alphabar))\right).$$
The identification of the right-hand side of the formula above with the definition of the rational function $B_{n, a_1, \dots, a_n}(v)(w, \plgt)$ (see Definition \ref{definitiongeometricbernoulli}) gives the conclusion:
$$\prod_{j= 1}^{n} G_{n-2, ((a_k)_{k \neq j})}(v)(w, x) ^{(-1)^{j+1}} =\exp(2i\pi B_{n, a_1, \dots, a_n}(v)(w, x)).$$

2. Consider $g \in \slnz{n}$. Remember that the actions of $\slnz{n}$ on $\homlc$ and $L$ satisfy $(g\cdot a) (g\cdot \alpha) = a(\alpha)$ for any $(a, \alpha) \in \homlc \times L$. The action on $L$ further extends to $V = L \otimes \qq$ and passes to the quotient $V/L$. Thus, in the construction of $G_{n-2, a_1, \dots, a_{n-1}}$ or $B_{n, a_1, \dots, a_n}$ replacing $\abar$ with $g \cdot \abar$ and $v$ with $g \cdot v$ replaces $\alphabar$ with $g\cdot \alphabar$, $\gamma$ with $g\cdot \gamma$ and $F = F(\abar, \alphabar, v)$ with $F_g = F(g\cdot\abar, g\cdot\alphabar, g\cdot v) = g \cdot F$ while $\epsilon$ is left unchanged. Therefore we may write:
$$G_{n-2, g\cdot\abar}(g \cdot v)(w, g\cdot \plgt) = \prod_{\delta \in F_g/\zz g\cdot\gamma} G_{n-2}\left(\frac{w + (g\cdot\plgt)(\delta)}{(g\cdot \plgt)(g\cdot\gamma)}, \frac{1}{(g\cdot \plgt)(g\cdot\gamma)} (g\cdot \plgt)(g\cdot\alphabar)\right). $$
Then, identifying $F_g = g\cdot F$ and putting $\delta' = g\cdot\delta$ gives
$$G_{n-2, g\cdot\abar}(g\cdot v)(w, g\cdot \plgt) = \prod_{\delta' \in F/\zz\gamma} G_{n-2}\left(\frac{w + \plgt(\delta')}{\plgt(\gamma)}, \frac{1}{\plgt(\gamma)} \plgt(\alphabar) \right) $$
which gives the conclusion:
$$G_{n-2, g\cdot\abar}(g\cdot v)(w, g\cdot \plgt) = G_{n-2, \abar}(v)(w, \plgt)$$
As for $B_{n, g\cdot a_1, \dots, g\cdot a_n}(g \cdot v)$, the set $F = F(\abar, \alphabar, v)$ is once again replaced by $g \cdot F$ so that:
$$B_{n, g\cdot a_1, \dots, g\cdot a_n}(g\cdot v)(w, g\cdot \plgt) = \frac{\epsilon}{n!} \sum_{\delta \in g\cdot F} B_{n,n}(w + (g\cdot \plgt)(\delta), (g\cdot \plgt)(g\cdot\alphabar)).$$
Put once again $\delta' = g \cdot \delta$, which gives:
$$B_{n, g\cdot a_1, \dots, g\cdot a_n}(g \cdot v)(w, g\cdot \plgt) = \frac{\epsilon}{n!} \sum_{\delta' \in F} B_{n,n}(w + \plgt(\delta'), \plgt(\alphabar)).$$
Identify the right-hand side to conclude that:
$$B_{n, g\cdot a_1, \dots, g\cdot a_n}(g \cdot v)(w, g\cdot \plgt) = B_{n, a_1, \dots, a_n}(v)(w, \plgt).$$
\end{proofbis}

\textbf{Remark:} The set $\mathcal{F}(V/L \times \cc \times \homlc, \cc)$ is naturally endowed with an action of $\slnz{n}$ given by $(g \cdot f)(v)(w, \plgt) = f(g^{-1}\cdot v)(w, g^{-1}\cdot \plgt)$. The second part of Theorem \ref{theoremmodulargeometricgrintroduction} may then be restated as:
\begin{align*}
g \cdot G_{n-2, a_1, \dots, a_{n-1}} &= G_{n-2,g \cdot a_1, \dots, g \cdot a_{n-1}}  \\
g \cdot B_{n, a_1, \dots, a_n} &= B_{n, g\cdot a_1, \dots, g\cdot a_n}
\end{align*}

We also add that both functions $G_{n-2, \abar}$ and $B_{n, \abar}$ behave nicely under permutation of vectors, namely for any permutation $\sigma \in \goth{S}_{n-1}$, $G_{n-2, \sigma(\abar)} = G_{n-2, \abar}^{\signature(\sigma)}$ where $\signature(\sigma)$ is the signature of the permutation $\sigma$ and for any permutation $\sigma \in \goth{S}_{n}$, $B_{n, \sigma(\abar)} = \signature(\sigma)B_{n, \abar}$. The modular property \refp{modulargeometricgrintroduction} may be restated as a partial coboundary relation between two collections of functions. Indeed, fixing a non zero primitive linear form $a \in \Lambda$ as a base point we may define:
$$\psi_{n,a}  := \begin{cases}\slnz{n}^{n-2} & \to \mathcal{F}(V/L \times \cc \times \homlc, \cc) \\ (g_1, \dots, g_{n-2}) & \to \left((v, w, \plgt) \to G_{n-2, a, g_1\cdot a, \dots, (g_1\dots g_{n-2}) \cdot a}(v)(w, \plgt)\right)\end{cases}$$
$$\phi_{n,a}  := \begin{cases} \slnz{n}^{n-1} & \to \mathcal{F}(V/L, \qq[w](\plgt)) \\ (g_1, \dots, g_{n-1}) & \to  B_{n, a, g_1\cdot a, (g_1g_2)\cdot a, \dots, (g_1\dots g_{n-1})\cdot a}(v)(w, \plgt)\end{cases}$$
When the linear forms $a, g_1 \cdot a, (g_1g_2)\cdot a, \dots, (g_1\dots g_{n-1}) \cdot a$ are linearly independent, we may rewrite formula $\refp{modulargeometricgrintroduction}$ as a relation between the multiplicative coboundary of $\psi_{n,a}$ and $\phi_{n,a}$ on a dense open subset of $\cc \times \homlc$:
\begin{equation}\label{cocyclegr}
\cohomdx\psi_{n,a}(g_1, \dots, g_{n-1})(v)(w,\plgt) = \exp(2i\pi \phi_{n,a}(g_1, \dots, g_{n-1})(v)(w, \plgt))
\end{equation}

We end this section by proving that the geometric families of $G_{r, \abar}$ functions satisfy distribution relations which generalise Proposition \ref{distributiongr}.

\begin{general}{Proposition}\label{distributiongrgeom}
Consider an integer $N \geq 2$. Then the following distribution relations hold:
\begin{equation}\label{formuladistributiongrgeom}\prod_{Nv' \equiv v \mod L} G_{r, \abar}(v')(w, \plgt) = G_{r, \abar}(v)(Nw,\plgt).
\end{equation}
\end{general}

\begin{proof}
We will rewrite the left-hand side of \refp{formuladistributiongrgeom} in order to use the results from Proposition \ref{distributiongr}. Indeed:
$$ \prod_{Nv' \equiv v \mod L} G_{r, \abar}(v')(w, \plgt) = \prod_{\beta \in L/NL} G_{r, \abar}\left(\frac{v + \beta}{N}\right)(w, \plgt).$$
Using the definition of $G_{r, \abar}$ we get:
$$ \prod_{Nv' \equiv v \mod L} G_{r, \abar}(v')(w, \plgt) = \prod_{\beta \in L/NL} \prod_{\delta' \in F(\abar, \alphabar, (v+\beta)/N)/\zz\gamma} G_r\left(\frac{w+\plgt(\delta')}{\plgt(\gamma)}, \frac{1}{\plgt(\gamma)}\plgt(\alphabar) \right).$$
Define $F_N = \disjointunion_{\beta \in L/NL} F(\abar, \alphabar, (v+\beta)/N)$. Using the definition of $F(\abar, \alphabar, \cdot)$ it is clear that:
\begin{align*}
F_N & = \disjointunion_{\beta \in L/NL}\left\{\delta' \in \frac{v+\beta}{N} + L \setseparator \forall\,1 \leq j \leq r+1,~0 \leq a_j(\delta) < a_j(\alpha_j) \right\} \\
F_N & = \left\{\delta' \in \frac{v+L}{N} \setseparator \forall\,1 \leq j \leq r+1,~0 \leq a_j(\delta) < a_j(\alpha_j) \right\}.
\end{align*}
Consider $\delta' \in F_N$. Then $N\delta'$ belongs to 
$$\left\{ \delta \in v + L \setseparator \forall\,1 \leq j \leq r+1,~0 \leq a_j(\delta) < Na_j(\alpha_j) \right\}$$ 
and there are unique integers $0 \leq k_1, \dots, k_{r+1} < N$ such that 
$$N\delta' - \sum_{j = 1}^{r+1} k_j \alpha_j \in \left\{ \delta \in v + L \setseparator \forall\,1 \leq j \leq r+1,~0 \leq a_j(\delta) < a_j(\alpha_j) \right\} = F(\abar, \alphabar, v).$$ 
Thus we get a bijection:
$$f := \begin{cases} F_N/\zz\gamma & \to \{0, 1, \dots, N-1\}^{r+2} \times F(\abar, \alphabar, v) \\ \delta' & \to ((k, k_1, \dots, k_{r+1}), \delta)\end{cases} $$
defined by $N\delta' \equiv k \gamma + \sum_{j = 1}^{r+1} k_j \alpha_j + \delta \mod \zz\gamma$. It follows that:
\begin{align*}
\prod_{Nv' \equiv v \mod L} G_{r, \abar}(v')(w, \plgt) & = \prod_{\delta' \in F_N/\zz\gamma}G_r\left(\frac{w+\plgt(\delta')}{\plgt(\gamma)}, \frac{1}{\plgt(\gamma)}\plgt(\alphabar) \right) \\
\prod_{Nv' \equiv v \mod L} G_{r, \abar}(v')(w, \plgt) & = \prod_{\delta \in F(\abar, \alphabar, v)/\zz\gamma} \prod_{k, k_1, \dots, k_{r+1} = 0}^{N-1} G_r\left(\frac{w+\plgt\left(\frac{\delta + k\gamma + k_1 \alpha_1 + \dots +k_{r+1}\alpha_{r+1}}{N}\right)}{\plgt(\gamma)}, \frac{1}{\plgt(\gamma)}\plgt(\alphabar) \right)
\end{align*}
Using the third relation from Proposition \ref{distributiongr} we get:
\begin{align*}
\prod_{Nv' \equiv v \mod L} G_{r, \abar}(v')(w, \plgt) & = \prod_{\delta \in F(\abar, \alphabar, v)/\zz\gamma} G_r\left(\frac{Nw+\plgt(\delta)}{\plgt(\gamma)}, \frac{1}{\plgt(\gamma)}\plgt(\alphabar) \right)\\
\prod_{Nv' \equiv v \mod L} G_{r, \abar}(v')(w, \plgt) & = G_{r, \abar}(v)(Nw, \plgt)
\end{align*}
which is the desired result.
\end{proof}

In an upcoming paper in this series, we use the $G_{r, \abar}$ functions to construct conjectural elliptic units above number fields with exactly one complex place which should behave as the Siegel units, and these distribution relations already show some of these similarities.

For the rest of this article we shift our focus from the $G_{r, \abar}$ functions to the collection of Bernoulli rational functions $B_{n, a_1, \dots, a_n}$ and the associated collection of $(n-1)$-cocycles $\phi_{n, a}$. In particular, we will show that formula \refp{cocyclegr} holds under less restrictive conditions on the $g_i$'s and that $\phi_{n,a}$ truly becomes a cocycle on specific subgroups of $\slnz{n}$.

\section{Cocycle properties for the collection of $B_{n, \abar}$ functions}\label{sectioncones}

The goal of this section is to show that the additive cocycle relation:
\begin{equation}
\sum_{j = 0}^{r+2} (-1)^j B_{n, a_0, \dots, \omitvar{a_j}, \dots, a_n}(v)(w, \plgt) \in \zz
\end{equation}
which holds for linear forms $a_0, \dots, a_n$ in general position in a rank $n$ lattice $\Lambda$ as a consequence of Theorem \ref{theoremmodulargeometricgrintroduction} may be improved to a finer cocycle relation:
\begin{equation}\label{cocyclerelationgeombern}
\sum_{j = 0}^{r+2} (-1)^j B_{n, a_0, \dots, \omitvar{a_j}, \dots, a_n}(v)(w, \plgt) = 0
\end{equation}
which holds for a wider range of configurations of $a_0, \dots, a_n$ inside the rank $n$ lattice $\Lambda$. This will be achieved in Proposition \ref{propositionbernoullicocycle} as we show that this relation may be obtained as the specialisation of a cocycle relation for indicator functions of closed cones. 

\subsection{A cocycle relation for closed cones}\label{sectionconescocycle}

We now introduce several notations for cones in a $\qq$-vector space $V$ of finite dimension $n$.  These notations as well as the strategies for the proofs by induction on the dimension are inspired by \cite{RichardHill}, \cite{CDG}, and are typical of the theory of polyhedral cones. We focus on $\qq$-vector spaces for ease of presentation, but we wish to highlight that these results would hold for a vector space over $\rr$ or over any ordered field. For a more complete presentation on the theory of convex cones we refer to \cite{BarvinokConvexity}. 

Let $V$ be a $\qq$-vector space of dimension $n$. A convex cone in $V$ is any convex set $C$ satisfying $\forall x, y \in C,  x+y \in C$ and $\forall \lambda > 0, \lambda.x \in C$. A convex cone $C$ is said to be polyhedral if there are two (possibly empty) sets of vectors $v_1, \dots, v_p$ and $v'_1, \dots, v'_q$ in $V$ such that $C = \qq_{\geq 0} v_1 + \dots +\qq_{\geq 0} v_p + \qq_{> 0} v'_1 + \dots \qq_{> 0} v'_q$. The vectors $v_i$ and $v'_j$ are called generators of the cone $C$. Note that by convention $\{0\}$ is a polyhedral cone with empty set of generators.  Let $\spancones$ be the $\qq$-algebra of $\qq$-valued functions on $V$ generated by the indicator functions of polyhedral cones. Denote by $\spanwedges$ the subspace of $\spancones$ generated by the indicator functions of the closed polyhedral cones containing a line $\qq v$ for some non zero vector $v \in V$. Note that $\spanwedges$ is not stable under multiplication. Most statements in the theory of cones may be proved by induction on the dimension and rest on the fact that when $V'$ is a subspace of $V$ there is a natural inclusion:
\begin{equation}\label{inclusionspancones}
\begin{cases} 
\prespancones{V'} & \to ~\spancones \\
f &\to ~\tilde{f} : v \to \begin{cases} f(v) & \text{ if } v \in V' \\ 0 & \text{ otherwise}\end{cases}
\end{cases}
\end{equation}
and this inclusion sends $\prespanwedges{V'}$ to $\spanwedges$. 

We now focus on closed cones. For any $v_1, \dots, v_m \in V$, we denote by $\cclosed(v_1, \dots, v_m)$ the indicator function of the closed polyhedral cone $\qq_{\geq 0} v_1+ \dots+ \qq_{\geq 0} v_m$. In this article, we also use a dual representation for the closed polyhedral cones.
Namely, for linear forms $a_1, \dots, a_m \in \dualsimple{V}$ the set 
$$ \cap_{i = 1}^m \{ v \in V \setseparator a_i(v) \geq 0\}$$
is a closed polyhedral cone in $V$ and we denote by $\cdual(a_1, \dots, a_m)$ its indicator function. Note that by convention we may put $\cclosed(\emptyset) = \dirac$ where $\dirac$ is the Dirac function at $0$ whereas $\cdual(\emptyset) = 1$ is the indicator function of $V$. This dual representation of a cone already allows us to shed a different light on the notion of positive dual family (see Definition \ref{definitionpositivedualfamily}). Indeed, fix $B = [e_1, \dots, e_n]$ a basis of $V$ and denote by $L$ the $\zz$-lattice $\oplus_{j = 1}^n \zz e_j$. It is clear that generators of cones may be rescaled, so that any polyhedral cone in $V$ admits a set of generators which lie inside $L$. It is also true that any closed polyhedral cone in $V$ admits a dual representation with linear forms $a_1, \dots, a_m \in \Lambda = \homlz$ for some integer $m \geq 0$. In the specific case where $m = n$ and $a_1, \dots, a_n \in \Lambda$ are linearly independent, the primitive positive dual family $\alpha_1, \dots, \alpha_n \in L$ to $a_1, \dots, a_n$ satisfies $\cdual(a_1, \dots, a_n) = \cclosed(\alpha_1, \dots, \alpha_n)$. 
Lastly, we say that a family $(v_1, \dots, v_m) \in V^m$ is in general position in $V$ if any of its subfamilies of size at most $n$ is free. 

In [\hspace{1sp}\cite{CDG}, \symbolparagraph1] Charollois, Dasgupta and Greenberg describe a cocycle relation for indicator functions of polyhedral cones using previous work by Hill \cite{RichardHill} on open polyhedral cones. Using the Solomon-Hu pairing (see \cite{SolomonHu}) they construct Shintani $(n-1)$-cocycles for $\slnz{n}$ with values in some spaces of rational functions which are cohomologous to Sczech cocycles \cite{Sczech}. Our goal is to describe another cocycle relation for indicator functions $\cdual(a_1, \dots, a_m)$ of closed polyhedral cones in a dual setting and show how formula \refp{cocyclerelationgeombern} may be deduced via the Solomon-Hu pairing. 
We start by stating some basic properties which will be very useful.

\begin{general}{Lemma}\label{lemmapropertiescones}
Let $V$ be a $\qq$-vector space of dimension $n$ and set $\dualsimple{V}$.
\begin{itemize}
\item[$\mathrm{(i)}$] If $a_1, \dots, a_m \in \dualsimple{V}$ do not generate $\dualsimple{V}$ (in particular if $m < n$) then $\cdual(a_1, \dots, a_m) \in \spanwedges$.
\item[$\mathrm{(ii)}$] For any $a_1, \dots, a_m \in \dualsimple{V}$, $\cdual(a_1, \dots, a_m) = \prod_{j = 1}^m \cdual(a_j)$.
\item[$\mathrm{(iii)}$] If $a, a' \in \dualsimple{V}$ satisfy $a = \lambda a'$ for some $\lambda > 0$ then $\cdual(a) = \cdual(a')$.
\item[$\mathrm{(iv)}$] For any $a \in \dualsimple{V}$, $\cdual(a, -a)$ is the indicator function of $\ker(a)$ and $\cdual(a, -a) + 1 = \cdual(a) + \cdual(-a)$.
\item[$\mathrm{(v)}$] If $a_1, \dots, a_m \in \dualsimple{V}$ generate $\dualsimple{V}$ and if there are positive coefficients $\lambda_1, \dots, \lambda_m$ such that $\sum_{j = 1}^{m} \lambda_ja_j = 0$, then $\cdual(a_1, \dots, a_m) = \dirac$ is the Dirac function at $0$.
\item[$\mathrm{(vi)}$] If $a_1, \dots, a_{m+1} \in \dualsimple{V}$ satisfy $\sum_{j = 1}^{m+1} \lambda_j a_j = 0$ with $\lambda_{m+1} < 0$ and $\lambda_j \geq 0$ for $1 \leq j \leq m$ then $\cdual(a_1, \dots, a_{m+1}) = \cdual(a_1, \dots, a_m)$.
\end{itemize}
\end{general}

\begin{proof}
(i) If $a_1, \dots, a_m$ do not generate $\dualsimple{V}$ then $\cap_{j = 1}^m \ker(a_j) \neq \{0\}$. For any $v \in \cap_{j = 1}^m \ker(a_j) - \{0\}$, and any $\lambda \in \qq$, it is clear that $\cdual(a_1, \dots, a_m)(\lambda.v) = 1$. Therefore the cone described by $\cdual(a_1, \dots, a_m)$ contains the line $\qq v$ and $\cdual(a_1, \dots, a_m) \in \spanwedges$. \smallskip

(ii) The cone described by $\cdual(a_1, \dots, a_m)$ is naturally defined as the intersection of the cones described by $\cdual(a_1), \dots, \cdual(a_m)$, therefore $\cdual(a_1, \dots, a_m) = \prod_{j = 1}^m \cdual(a_j)$. \smallskip

(iii) If $a = \lambda a'$ with $\lambda > 0$ then for any $v \in V$, $a(v) \geq 0 \Leftrightarrow a'(v) \geq 0$, which gives $\cdual(a) = \cdual(a')$. \smallskip

(iv) For $a \in \dualsimple{V}$ and $v \in V$, $\cdual(a, -a)(v) = 1$ if and only if $a(v) \geq 0$ and $-a(v) \geq 0$ i.e. if and only if $a(v) = 0$. The equality $\cdual(a, -a)(v) + 1 = \cdual(a)(v) + \cdual(-a)(v)$ is easily computed in all three cases $a(v) >0$, $a(v) = 0$ and $a(v) < 0$. \smallskip

(v) Suppose that $a_1, \dots, a_m$ generate $\dualsimple{V}$ and that there are coefficients $\lambda_j > 0$ such that $\sum_{j = 1}^m \lambda_j a_j = 0$. Suppose that there is a non-zero vector $v \in V$ such that $\cdual(a_1, \dots, a_m) = 1$. Then $a_j(v) \geq 0$ for $1 \leq j \leq m$. As $a_1, \dots, a_m$ generate $\dualsimple{V}$ and $v \neq 0$, there is an index $1 \leq l \leq m$ such that $a_l(v) > 0$. The sum $\sum_{j = 1}^m \lambda_ja_j(v)$ is equal to $0$ and contains only non-negative terms, therefore all terms must be zero, which contradicts $\lambda_l a_l(v) > 0$. Then we only check that $\cdual(a_1, \dots, a_m)(0) = 1$ and conclude that $\cdual(a_1, \dots, a_m) = \dirac$ is the Dirac function at $0$.

(vi) Suppose $\sum_{j = 1}^{m+1} \lambda_j a_j = 0$ with $\lambda_j \geq 0$ for $1 \leq j \leq m$ and $\lambda_{m +1} < 0$. Then $a_1(v) \geq 0, \dots, a_m(v) \geq 0 \Rightarrow a_{m+1}(v) \geq 0$ which gives 
$$\cap_{i = 1}^m \{ v \in V \setseparator a_i(v) \geq 0\} = \cap_{i = 1}^{m+1} \{ v \in V \setseparator a_i(v) \geq 0\}$$ 
and therefore $\cdual(a_1, \dots, a_{m+1}) = \cdual(a_1, \dots, a_m)$.
\end{proof}

We now give a crucial definition which will be used in the proof of Theorem  \ref{theoremcocyclecones}.

\begin{general}{Definition}\label{definitionrelation}
Let $V$ be a $\qq$-vector space and let $a_0, \dots, a_m$ be $m+1$ non-zero linear forms on $V$ such that $\rank(a_0, \dots, a_m) = m$. There is a unique linear combination $\sum_{j = 0}^m \lambda_j a_j = 0$ with coefficients $\lambda_j \in \qq$ satisfying:
\begin{itemize}
\item $(\lambda_0, \dots, \lambda_m) \neq (0, \dots, 0)$
\item $\cardinalshort{\{0 \leq j \leq m \setseparator \lambda_j < 0\}} \leq \cardinalshort{\{0 \leq j \leq m \setseparator \lambda_j > 0\}}$
\item if $\cardinalshort{\{0 \leq j \leq m \setseparator \lambda_j < 0\}} = \cardinalshort{\{0 \leq j \leq m \setseparator \lambda_j > 0\}}$, the first non-zero coefficient $\lambda_l$ is negative
\item the first non-zero coefficient $\lambda_l$ has absolute value $1$.
\end{itemize}
This linear combination will be refered to as the standard non-trivial relation among $a_0, \dots, a_m$. In this situation, we define $k^{-}(a_0, \dots, a_m) = \cardinalshort{\{0 \leq j \leq m \setseparator \lambda_j < 0\}}$. 
\end{general}

In particular, when $V$ has finite dimension $n$, Definition \ref{definitionrelation} applies to any family $a_0, \dots, a_n \in \dualsimple{V}$ which generates $\dualsimple{V}$. In this situation, we may also define $k^{0}(a_0, \dots, a_m)$ (resp. $k^{+}(a_0, \dots, a_m)$) the number of coefficients $\lambda_j = 0$ (resp. $\lambda_j > 0$) in the relation, but these won't be much needed. Before giving the proof of Theorem \ref{theoremcocyclecones}, we discuss the configurations of linear forms to which it applies in the light of Definition \ref{definitionrelation}. Indeed, we have carefully avoided some configurations of the linear forms $a_j$ which we will refer to as ``bad position'' or \refp{badpositioncondition} for short. Namely, in the vector space $\dualsimple{V}$ of dimension $n$, a set of non-zero linear forms $a_0, \dots, a_n$ is in bad position \refp{badpositioncondition} if:
\begin{equation}\label{badpositioncondition}\tag{BP}
\rank(a_0, \dots, a_n) = n \text{ and } k^{-}(a_0, \dots, a_n) = 0 \text{ and } k^0(a_0, \dots, a_n) > 0.
\end{equation}
Theorem \ref{theoremcocyclecones} applies to all other configurations of the linear forms
$a_0, \dots, a_n$ with a specific treatment when $\rank(a_0, \dots, a_n) = n$ and $k^{+}(a_0, \dots, a_n) = n+1$. The specific behaviour associated to this configuration regarding cocycle relations was already observed in \mbox{[\hspace{1sp}\cite{RichardHill}, Proposition 2]}.

We are now ready to prove Theorem \ref{theoremcocyclecones}, showing that the functions $\cdual(\cdot)$ satisfy some cocycle relations. These are inspired by the cocycle relations described by Hill [\hspace{1sp}\cite{RichardHill}, Proposition 2] for the indicator functions $\copen(v_1, \dots, v_m)$ of open cones $\rr_{> 0} v_1 + \dots + \rr_{> 0} v_m$ and by Charollois, Dasgupta and Greenberg [\hspace{1sp}\cite{CDG}, Theorem 1.1] for variants of the functions $\copen(\cdot)$ for which some boundaries are included depending on a so-called $Q$-perturbation process. We argue that Theorem \ref{theoremcocyclecones} is simpler as there is no need to select the boundary pieces to add or remove with such a process, yet the proof uses 
similar ideas as those developed in the proof of [\hspace{1sp}\cite{CDG}, Theorem 1.1]. 

\bigskip

\begin{proofbis}{Proof of Theorem \ref{theoremcocyclecones}}
Let $V$ be a $\qq$-vector space of dimension $n$. Consider $n+1$ linear forms $a_0, \dots, a_n \in \dualsimple{V}$ which generate $\dualsimple{V}$ and are not \refp{badpositioncondition}. Define $q(k)$ by $q(0) = 1$ and $q(k) = 0$ otherwise. We wish to prove that:
$$ \sum_{j = 0}^n \eps_j\cdual(a_0, \dots, \omitvar{a_j}, \dots, a_n) \equiv q(k^{-}(a_0, \dots, a_n)) \eps_0\dirac \mod \spanwedges $$
where $\eps_j = (-1)^j \signdet(a_0, \dots, \omitvar{a_j}, \dots, a_n)$ and $\dirac$ is the Dirac at $0$. The proof is split into three parts. We first treat the case where $a_0, \dots, a_n$ are in general position with $k^{-}(a_0, \dots, a_n) = 0$. Then, using this first result we treat the case where $a_0, \dots, a_n$ are in general position with $k^{-}(a_0, \dots, a_n) >0$ by double induction on the dimension $n$ and the value of $k^{-}(a_0, \dots, a_n)$. Finally, using this second result, we treat the case where $a_0, \dots, a_n$ are not in general position and $k^{-}(a_0, \dots, a_n) > 0$ by single induction on $n$. \bigskip

\underline{First case: $a_0, \dots, a_n$ are in general position and $k^{-}(a_0, \dots, a_n) = 0$:} \smallskip
\newline Let $\sum_{j = 0}^n \lambda_j a_j = 0$ be the standard non-trivial relation among $a_0, \dots, a_n$. Because the coefficients $\lambda_j$ are all positive, the signs $\eps_j$ are all equal to $\eps_0$. Indeed, if $j \geq 1$ then:
\begin{align*}
\eps_j &= (-1)^j \signdet(-\lambda_ja_j/\lambda_0, \dots, \omitvar{a_j}, \dots, a_{n}) \\
\eps_j &= -(-1)^j\sign(\lambda_j/\lambda_0)\signdet(a_j, a_1, \dots, \omitvar{a_j}, \dots, a_n) \\
\eps_i &= -(-1)^j(-1)^{j+1}\signdet(a_1, \dots, a_n)\\
\eps_j &= \eps_0
\end{align*}
Let us denote by $\eps$ the common sign of the $\eps_0 = \dots = \eps_n = \eps$. Consider now $f(a_0, \dots, a_n) = \prod_{j = 0}^n (\cdual(a_j)-1)$. Expanding the product gives:
\begin{align*}
f(a_0, \dots, a_n) &= \sum_{P \subset [|0, n|]} (-1)^{n+1-\cardinalshort{P}} \prod_{j \in P} \cdual(a_j)\\
f(a_0, \dots, a_n) &= \sum_{P \subset [|0, n|]} (-1)^{n+1-\cardinalshort{P}} \cdual(a_j, j \in P)
\end{align*}
The last line follows from lemma \ref{lemmapropertiescones} (ii). For any $P \subset [|0, n|]$ with $\cardinalshort{P} < n$, the function $\cdual(a_j,  j\in P)$ belongs to $\spanwedges$ (lemma \ref{lemmapropertiescones} (i)). Thus:
$$f(a_0, \dots, a_n) -\cdual(a_0, \dots, a_n) +\sum_{j = 0}^n \cdual(a_0, \dots, \omitvar{a_j}, \dots, a_n) \in \spanwedges$$
Since $\sum_{j = 0}^n \eps_j \cdual(a_0, \dots, \omitvar{a_j}, \dots, a_n) = \eps\sum_{j = 0}^n \cdual(a_0, \dots, \omitvar{a_j}, \dots, a_n)$ this gives:
$$\eps f(a_0, \dots, a_n) -\eps\cdual(a_0, \dots, a_n) +\sum_{j = 0}^n \eps_j\cdual(a_0, \dots, \omitvar{a_j}, \dots, a_n) \in \spanwedges$$
Now, since the coefficients $\lambda_i$ are all positive by assumption, $\cdual(a_0, \dots, a_n)$ is the dirac at $0$ by lemma \ref{lemmapropertiescones} (v). On the other hand:
$$f(a_0, \dots, a_n) = \prod_{j = 0}^n (\cdual(a_j) -1) = 0$$
The last equality holds because for all $v \in V$, $\sum_{j = 0}^n \lambda_j a_j(v) = 0$ so at least one of the $a_j(v)$ must be non-negative. Therefore we may conclude that $\sum_{j = 0}^n \eps_j \cdual(a_0, \dots, \omitvar{a_j}, \dots, a_n) \equiv \eps \dirac \mod \spanwedges$ as claimed. \bigskip

\underline{Second case: $a_0, \dots, a_n$ are in general position and $k^{-}(a_0, \dots, a_n) > 0$:} \smallskip
\newline the case $n = 1$ is immediate as in this case $\lambda_0 a_0 + \lambda_1 a_1 = 0$ with $\lambda_0 = -1$ and $\lambda_1 > 0$ which gives $\cdual(a_0) - \cdual(a_1) = \cdual(a_0) - \cdual(a_0) = 0 \in \spanwedges$ by lemma \ref{lemmapropertiescones} (iii). We are now ready to perform double induction on both $n$ and $k^{-}$. Suppose that the result holds for any family $a_0, \dots, a_{n-1}$ of linear forms on a $n-1$ dimensional $\qq$-vector space $V'$ generating $\dualsimple{V'}$ with $k^{-}(a_0, \dots, a_{n-1}) > 0$ and that it holds in dimension $n$ whenever $a_0, \dots, a_n$ are in general position in $\dualsimple{V}$ with $k^{-}(a_0, \dots, a_n) = k^{-}$. Suppose now that $a_0, \dots, a_n$ are in general position in $\dualsimple{V}$ with $k^{-}(a_0, \dots, a_n) = k^{-}+1$. We aim to prove that:
$$ \sum_{j = 0}^n \eps_j\cdual(a_0, \dots, \omitvar{a_j}, \dots, a_n) \in \spanwedges $$
where as before $\eps_j = (-1)^j \signdet(a_0, \dots, \omitvar{a_j}, \dots, a_n)$. Denote as before $\sum_{j = 0}^n \lambda_j a_j = 0$ the standard non-trivial relation among $a_0, \dots, a_n$. By assumption, $k^{-}(a_0, \dots, a_n) > 0$ so there is at least one index $l$ such that $\lambda_l < 0$. We fix any such index $l$ and we will show that the desired result may be deduced from the result for the families $a_0, \dots, a_{l-1}, -a_l, a_{l+1}, \dots, a_n$ in $\dualsimple{V}$ and $a_{0\at{\ker{a_l}}}, \dots, \omitvar{a_{l\at{\ker{a_l}}}}, \dots, a_{n\at{\ker{a_l}}}$ in $\dualsimple{(\ker{a_l})}$ as follows. For simplicity we use once again the auxiliary function defined by $q(k) = 1$ if $k = 0$ and $q(k) = 0$ otherwise. The family $(a'_0, \dots, a'_n) = (a_0, \dots, a_{l-1}, -a_l, a_{l+1}, \dots, a_n)$ is in general position in $\dualsimple{V}$ with $k^{-}(a'_0, \dots, a'_n) = k^{-}$. By induction hypothesis on $k^{-}$, fixing any index $m \neq l$:
\begin{equation}\label{eqstar}
\sum_{j = 0}^n \eps'_j\cdual(a'_0, \dots, \omitvar{a'_j}, \dots, a'_n) \equiv q(k^{-})\eps'_m \dirac \mod \spanwedges
\end{equation}
where $\eps'_j =(-1)^j\signdet(a'_0, \dots, \omitvar{a'_j}, \dots, a'_n)$. It is clear that $\eps'_j = -\eps_j$ when $j \neq l$ and $\eps'_l = \eps_l$ so that \refp{eqstar} reads:
\begin{equation}\label{eqstara}
\eps_l\cdual(a_0, \dots, \omitvar{a_l}, \dots, a_n) + \sum_{\substack{j = 0\\ j \neq l}}^n (-\eps_j)\cdual(a'_0, \dots, \omitvar{a'_j}, \dots, a'_n) \equiv -q(k^{-})\eps_m\dirac \mod \spanwedges.
\end{equation}
It is then sufficient to prove that:
$$ \left(\sum_{\substack{j = 0\\j \neq l}}^n\eps_j(\cdual(a_0, \dots, \omitvar{a_j}, \dots, a_n) +\cdual(a'_0, \dots, \omitvar{a'_j}, \dots, a'_n))\right) \equiv q(k^{-})\eps_m\dirac \mod \spanwedges$$
Using lemma \ref{lemmapropertiescones} (ii) then (iv), we obtain:
\begin{align*}
&\sum_{\substack{j = 0\\j \neq l}}^n\eps_j(\cdual(a_0, \dots, \omitvar{a_j}, \dots, a_n) +\cdual(a'_0, \dots, \omitvar{a'_j}, \dots, a'_n))\\
&=\sum_{\substack{j = 0\\j \neq l}}^n\eps_j(\cdual(a_0, \dots, \omitvar{a_j}, \dots, a_n) +\cdual(a_0, \dots, -a_l, \dots, \omitvar{a_j}, \dots, a_n))\\
&=\sum_{\substack{j = 0\\j \neq l}}^n\eps_j\cdual(a_0, \dots, \omitvar{a_l}, \dots, \omitvar{a_j}, \dots, a_n)(\cdual(a_l) + \cdual(-a_l))\\
&=\sum_{\substack{j = 0\\j \neq l}}^n\eps_j\cdual(a_0, \dots, \omitvar{a_l}, \dots, \omitvar{a_j}, \dots, a_n)(1+\ker(a_l))
\end{align*}
where $\ker(a_l)$ is the indicator function of the kernel of $a_l$. For any $j \neq l$, the function $\cdual(a_0, \dots, \omitvar{a_l}, \dots, \omitvar{a_j}, \dots, a_n)$ lies in $\spanwedges$ by lemma \ref{lemmapropertiescones} (i), thus we only need to prove that:
$$\sum_{\substack{j = 0\\j \neq l}}^n\eps_j\cdual(a_0, \dots, \omitvar{a_l}, \dots, \omitvar{a_j}, \dots, a_n)\ker(a_l) \equiv q(k^{-})\eps_m\dirac \mod \spanwedges$$
This step of the proof makes use of the induction hypothesis on $n$. Denote $a''_0, \dots, a''_n$ the restrictions of $a_0, \dots, a_n$ on $\ker(a_l)$. Then the natural inclusion described by \refp{inclusionspancones} gives:
$$\sum_{\substack{j = 0\\j \neq l}}^n\eps_j\cdual(a_0, \dots, \omitvar{a_l}, \dots, \omitvar{a_j}, \dots, a_n)\ker(a_l) = \sum_{\substack{j = 0\\j \neq l}}^n\eps_j\cdual_{\ker(a_l)}(a''_0, \dots, \omitvar{a''_l}, \dots, \omitvar{a''_j}, \dots, a''_n)$$
The elements $a''_0, \dots, \omitvar{a''_l}, \dots, a''_n$ are in general position in the dual $\dualsimple{\ker(a_l)}$ of $\ker(a_l)$. This space may be oriented using the form:
$$\mathcal{O}_l(x_0, \dots, \omitvar{x_l}, \dots, x_n) = \det(\tilde{x}_0, \dots, \tilde{x}_{l-1}, a_l, \tilde{x}_{l+1}, \dots, \tilde{x}_n)$$
for any lifts $\tilde{x}_j$ of the $x_j$ to $\dualsimple{V}$ satisfying $\tilde{x}_{j\at{\ker{a_l}}} = x_j$. Lastly, if $\sum_{j = 0}^n \lambda_j a_j = 0$ is the standard non-trivial relation among $a_0, \dots, a_n$ in $\dualsimple{V}$ then $\sum_{j \neq l} \lambda_j a''_j = 0$ is the standard non-trivial relation among $a''_0, \dots, \omitvar{a''_l}, \dots, a''_n$ in $\dualsimple{\ker(a_l)}$ and $k^{-}(a''_0, \dots, \omitvar{a''_l}, \dots, a''_n) = k^{-}(a_0, \dots, a_l, \dots, a_n)-1 = k^{-}$.
Therefore, using the induction hypothesis on $n$, we obtain:
$$\sum_{\substack{j = 0\\j \neq l}}^n\eps_j\cdual_{\ker(a_l)}(a''_0, \dots, \omitvar{a''_l}, \dots, \omitvar{a''_j}, \dots, a''_n) \equiv q(k^{-})\eps_m\dirac_{\ker(a_l)} \mod \mathcal{L}(\ker(a_l)) $$
where $\dirac_{\ker(a_l)}$ is the Dirac function at $0$ on $\ker(a_l)$. This last expression can be lifted back to $V$ using the inclusion \refp{inclusionspancones} as:
\begin{equation}\label{eqstarb}
\sum_{\substack{j = 0\\j \neq l}}^n\eps_j\cdual(a_0, \dots, \omitvar{a_l}, \dots, \omitvar{a_j}, \dots, a_n)\ker(a_l) \equiv q(k^{-})\eps_m\dirac \mod \spanwedges
\end{equation}
where it is clear that the inclusion map sends $\dirac_{\ker(a_l)}$ to $\dirac$.
We now piece \refp{eqstara} and \refp{eqstarb} together and find:
\begin{align*}
\sum_{j = 0}^n & (-1)^j  \signdet(a_0, \dots, \omitvar{a_j}, \dots, a_n)\cdual(a_0, \dots, \omitvar{a_j}, \dots, a_n) \\
& =\left(\sum_{j = 0}^n (-1)^j \signdet(a'_0, \dots, \omitvar{a'_j}, \dots, a'_n)\cdual(a'_0, \dots, \omitvar{a'_j}, \dots, a'_n) + q(k^{-})\eps_m\dirac\right)\\
& ~~+\left(\sum_{\substack{j = 0\\j \neq l}}^n\eps_j\cdual(a_0, \dots, \omitvar{a_l}, \dots, \omitvar{a_j}, \dots, a_n)\ker(a_l) -q(k^{-})\eps_m\dirac\right)
\end{align*}
Both terms in the right-hand side belong to $\spanwedges$, therefore the left hand side does too, as claimed. \bigskip

\underline{Third case: $a_0, \dots, a_n$ are not in general position and $k^{-}(a_0, \dots, a_n) > 0$:} \smallskip
\newline We prove this case by induction over $n \geq 2$ because it can't occur for $n = 1$. For $n = 2$ it may only happen when $a_j = \lambda a_l$ for some $j \neq l$ and some $\lambda \in \qq_{>0}$. Without loss of generality, assume $a_0 = \lambda a_1$. In that case, it follows from \ref{lemmapropertiescones} (ii) and (iii) that:
$$ \signdet(a_1, a_2)\cdual(a_1, a_2) - \signdet(a_0, a_2) \cdual(a_0, a_2) + 0 = 0 $$
which proves the case $n = 2$. Let us now assume the result holds in dimension $n-1$, where $n \geq 3$. Consider $V$ a $\qq$-vector space of dimension $n$ and suppose that $a_0, \dots, a_n$ are non-zero linear forms on $V$ which generate $\dualsimple{V}$. Assume that $a_0, \dots, a_n$ are not in general position in $\dualsimple{V}$ and that $k^{-}(a_0, \dots, a_n) > 0$. Without loss of generality we may assume that $a_n$ lies outside the span of $a_0, \dots, a_{n-1}$. Choose a vector $v_n \in \cap_{j = 0}^{n-1} \ker(a_j)$ such that $a_n(v_n) > 0$. This is possible because $a_0, \dots, a_{n-1}$ do not generate $\dualsimple{V}$. Let us now use the projection $ \pi : V \to \ker(a_n)$ with kernel $\qq v_n$ corresponding to the decomposition $V = \ker(a_n) \oplus \qq v_n$. For $0 \leq j \leq n$, denote $a'_j$ the restriction of $a_j$ to $\ker(a_n)$. Write as before $\eps_j = (-1)^j\signdet(a_0, \dots, \omitvar{a_j}, \dots, a_n)$. As $a_0, \dots, a_{n-1}$ are linearly dependent, $\eps_n = 0$, and we may compute:
$$ \sum_{j = 0}^n \eps_j\cdual(a_0, \dots, \omitvar{a_j}, \dots, a_n) = \sum_{j = 0}^{n-1} \eps_j\cdual(a_n)\cdual(a'_0, \dots, \omitvar{a'_j}, \dots, a'_{n-1})(\pi(\cdot)) $$
Let $\sum_{j = 0}^n \lambda_j a_j = 0$ be the standard non-trivial relation among $a_0, \dots, a_n$. The assumption that $a_n$ lies outside the span of $a_0, \dots, a_{n-1}$ is equivalent to $\lambda_n = 0$. Therefore, the standard non-trivial relation among $a'_0, \dots, a'_{n-1}$ is $\sum_{j = 0}^{n-1} \lambda_j a'_j = 0$ and $k^{-}(a'_0, \dots, a'_{n-1}) = k^{-}(a_0, \dots, a_n) > 0$. If the linear forms $a'_0, \dots, a'_{n-1}$ are in general position in $\dualsimple{(\ker{a_n})}$, we may use the result proven in the second case, and otherwise use the induction hypothesis as $a'_0, \dots, a'_{n-1}$ are not in general position in $\dualsimple{\ker(a_n)}$. In both cases, we obtain that:
$$ \sum_{j = 0}^{n-1} \eps_j\cdual(a'_0, \dots, \omitvar{a'_j}, \dots, a'_{n-1}) \in \mathcal{L}(V')$$
and remark that $\forall\,f \in \mathcal{L}(\ker(a_n))$, $\cdual(a_n)f(\pi(\cdot)) \in \spanwedges$ for any linear projection $\pi : V \to \ker(a_n)$. Indeed, if $f$ is the characteristic function of a cone containing a line $\qq v$ in $\ker(a_n)$, then $f(\pi(\cdot))$ is the characteristic function of a cone containing the line $\qq v$ in $V$. For $\lambda \in \qq$, $a_n(\lambda.v) = 0$ therefore $\cdual(a_n)(\lambda.v)f(\pi(\lambda.v)) = f(\lambda.v) = 1$ and $\cdual(a_n)f(\pi(\cdot))$ is the characteristic function of a cone containing the line $\qq v$. Thus:
$$\sum_{j = 0}^{n-1} \eps_j\cdual(a_n)\cdual(a'_0, \dots, \omitvar{a'_j}, \dots, a'_{n-1})(\pi(\cdot)) \in \spanwedges$$
 from which we conclude that $ \sum_{j = 0}^n \eps_j\cdual(a_0, \dots, \omitvar{a_j}, \dots, a_n) \in \spanwedges$ as claimed.
\end{proofbis}

\noindent\textbf{Remark:} In the case where $a_0, \dots, a_n$ do not span $\dualsimple{V}$ it is clearly also true that:
$$  \sum_{j = 0}^n (-1)^j\signdet(a_0, \dots, \omitvar{a_j}, \dots, a_n)\cdual(a_0, \dots, \omitvar{a_j}, \dots, a_n) =0$$
as each term in the left hand side is zero.

Unfortunately, the case where $a_0, \dots, a_n$ are not in general position with $k^{-}(a_0, \dots, a_n) = 0$ (this corresponds to the ``bad position'' condition \refp{badpositioncondition}) already fails in dimension $2$. Indeed, if $a_0, a_1, a_2$ satisfy $\rank(a_0, a_1, a_2) = 2$ and $\lambda_0a_0 + \lambda_1a_1 = 0$ with $\lambda_0 > 0$ and $\lambda_1 > 0$ then:
\begin{align*}
&\signdet(a_1, a_2)\cdual(a_1, a_2) - \signdet(a_0, a_2) \cdual(a_0, a_2) + \signdet(a_0, a_1) \cdual(a_0, a_1) \\
& =\signdet(-a_0, a_2)\cdual(-a_0, a_2) - \signdet(a_0, a_2) \cdual(a_0, a_2) + 0\\
& = -\signdet(a_0, a_2)\cdual(a_2)(\cdual(-a_0) + \cdual(a_0))\\
& = -\signdet(a_0, a_2)\cdual(a_2)(1 + \ker(a_0)) \not\in \spanwedges \oplus \qq \dirac
\end{align*}
Indeed, the function $\cdual(a_2) \in \spanwedges$ by lemma \ref{lemmapropertiescones} (i), however, the function $\cdual(a_2)\ker(a_0)$ is the indicator function of the cone $\qq_{\geq 0} v_2$ where $v_2$ is the vector defined by $a_0(v_2) = 0$ and $a_2(v_2) = 1$. 

The first part of the proof of Theorem \ref{theoremcocyclecones} may be slightly adjusted to prove that when $k^{-}(a_0, \dots, a_n) = 0$ and $a_0, \dots, a_n$ are not in general position, the following holds:
\begin{multline*}
$$ \sum_{i = 0}^n (-1)^i\signdet(a_0, \dots, \omitvar{a_i}, \dots, a_n)\cdual(a_0, \dots, \omitvar{a_i}, \dots, a_n) \equiv \\
\pm \prod_{\lambda_i = 0} \cdual(a_i) \indicator{\cap_{\lambda_i > 0} \ker(a_i)} \mod \spanwedges $$
\end{multline*}
It is not hard to see that the right-hand side of this formula doesn't belong to $\spanwedges \oplus \qq \dirac$. We will now introduce the generating functions of cones and relate this cocycle for indicator functions of cones to our geometric Bernoulli cocycle.

\subsection{Generating functions of cones and the $B_{n,\abar}$ functions}\label{sectionconesgenerating}

In this section we briefly recall some results on the generating functions associated to cones in $\rr^n$. Let $V_{\rr}$ be a $\rr$-vector space of finite dimension $n$ and fix an isomorphism $V_{\rr} \simeq \rr^n$ which defines $V_{\qq} = V \simeq \qq^n$ and $L = V_{\zz} \simeq \zz^n$. It is well-known that a generating function may be associated to any rational polyhedral cone in $V \simeq \rr^n$ which does not contain any line by the formula:
$$ g(C, v)(y) = \sum_{\delta \in C \cap (v+L)} y^\delta$$
where $y^\delta = y_1^{\delta_1}y_2^{\delta_2}\dots y_n^{\delta_n}$ and $v \in V_{\rr}$. This is well-defined for $y$ in a certain open subset of $\cc^n$ and it is possible to extend this function by analytic continuation. Indeed, if $C = \rr_{\geq 0} \alpha_1 \dots, \rr_{\geq 0} \alpha_n$ with primitive vectors $\alpha_1, \dots, \alpha_m \in L$ then the generating series associated to $C$ is in fact defined by:
$$ g(C, v)(y) = \frac{\sum_{\delta \in P \cap (v+L)}y^{\delta}}{(1-y^{\alpha_1})\dots(1-y^{\alpha_m})}$$
where $P = P(\alpha_1, \dots, \alpha_m) = \{ \sum_{i = 1}^m \mu_i \alpha_i \setseparator 0 \leq \mu_i < 1, \forall i\}$. This formula is given by the standard decomposition of a rational polyhedral cone:
\begin{equation}\label{decompositioncone}
C \cap (v + L) = \sqcup_{\delta \in P \cap (v+ L)} \delta + \zz_{\geq 0} \alpha_1 + \dots + \zz_{\geq 0} \alpha_m
\end{equation}
We highlight that Definitions \ref{propositiondefgeometricgr} and \ref{definitiongeometricbernoulli} used this decomposition implicitly with the set $F(\abar, \alphabar, v) = P(\alphabar) \cap (v +L)$. The key result regarding generating series associated to cones is the fact that the function $g$ may be extended to the subspace $\spanqconesrr$ of $\spanconesrr$ spanned by the indicator functions of rational polyhedral cones by linearity, as was proven independently by Khovanskii and Pukhlikov, and by Lawrence. A polyhedral cone in $V_{\rr}$ is said to be rational if it admits a set of generators inside $V_{\qq}$, in which case it also admits a set of generators inside $L = V_{\zz}$. The identification $\spanqconesrr \simeq \spancones \otimes_{\qq} \rr$ where $V = V_{\qq}$ thus shows that the generating series function $g$ may be extended to $\spancones$. In particular, for any $f \in \spanqwedgesrr \simeq \spanwedges \otimes_{\qq} \rr, g(f,\cdot) = 0$. The Solomon-Hu \cite{SolomonHu} pairing is then defined for $(f, v, \plgt) \in \spancones/\spanwedges \times V/L \times \homlc$ by $h(f, v)(\plgt) = g(f, v)(e^{\plgt})$. We use this pairing to express the geometric Bernoulli rational functions via a coefficient extraction.

\begin{general}{Definition}\label{definitionhzero}
The map
$$ h_0 := \begin{cases}\spancones & \to \mathcal{F}(V/L, \qq[w](\plgt)) \\ f & \to \coefficient{e^{wt}h(f, v)(t.\plgt)}{t}{0} \end{cases}$$
is $\qq$-linear and vanishes on $\spanwedges$.
\end{general}

It is also clear that $h_0(\dirac, v)$ is the constant function equal to $1$ if $v \in L$ and $0$ otherwise. We shall now use this function to express the geometric Bernoulli rational functions (see Definition \ref{definitiongeometricbernoulli}) associated to non-zero integral linear forms $a_1, \dots, a_n \in \Lambda= \homlz$.

\begin{general}{Lemma}\label{lemmabernoullicones}
For non-zero linear forms $a_1, \dots, a_n \in \Lambda$  in the rank $n$ oriented lattice $\Lambda$, for a vector $v \in V/L$ and for $(w, x)$ in a dense open subset of $\cc \times \homlc$:
$$B_{n, a_1, \dots, a_n}(v)(w, \plgt) = h_0(\signdet(a_1, \dots, a_n)\cdual(a_1, \dots, a_n), v)(w, \plgt)$$
\end{general}

\begin{proof}
If $a_1, \dots, a_n$ are linearly dependent, then both sides are zero so the equality holds. If $a_1, \dots, a_n$ are linearly independent, the proof follows from the definition of $B_{n, a_1, \dots, a_n}$: let $\alpha_1, \dots, \alpha_n$ be the primitive positive dual family to $a_1, \dots, a_n$ in $L \simeq \zz^n$. Denote $\epsilon = \signdet(a_1, \dots, a_n)$. Then identifying $F(\abar, \alphabar, v) = P(\alphabar) \cap (v + L)$ we get:
\begin{align*}
B_{n, a_1, \dots, a_n}(w, \plgt) &= \coefficient{(-1)^n \epsilon \sum_{\delta \in F(\abar, \alphabar, v)} \frac{e^{wt}e^{\plgt(\delta)t}}{\prod_{j = 1}^n(e^{\plgt(\alpha_j)t}-1)}}{t}{0}\\
B_{n, a_1, \dots, a_n}(w, \plgt) &= \coefficient{\epsilon \sum_{\delta \in P(\alphabar) \cap (v+L)} \frac{e^{wt}e^{\plgt(\delta)t}}{\prod_{j = 1}^n(1-e^{\plgt(\alpha_j)t})}}{t}{0}\\
B_{n, a_1, \dots, a_n}(w, \plgt) &= \coefficient{\epsilon e^{wt}h(\cclosed(\alpha_1, \dots, \alpha_{n}), v)(t.\plgt)}{t}{0}\\
B_{n, a_1, \dots, a_n}(w, \plgt) &=\coefficient{\epsilon e^{wt}h(\cdual(a_0, \dots, a_{r+1}), v)(t.\plgt)}{t}{0}
\end{align*}
where we use that $\cclosed(\alpha_1, \dots, \alpha_n) = \cdual(a_1, \dots, a_n)$ by definition of $\alpha_1, \dots, \alpha_n$. The identification of the right-hand side with $h_0(\epsilon\cdual(a_1, \dots, a_n), v)(w, \plgt)$ proves the claim.
\end{proof}

Using this lemma, we may finally describe the cocycle relations satisfied by the Bernoulli rational functions $B_{n, a_1, \dots, a_n}$.

\begin{general}{Proposition}\label{propositionbernoullicocycle}
Let $a_0, \dots, a_n$ be $n+1$ non-zero linear forms in $\Lambda$. Fix $v \in V/L$.
\begin{itemize}
\item Suppose $\rank(a_0, \dots, a_n) \leq n-1$, or $\rank(a_0, \dots, a_n) = n$ and $k^{-}(a_0, \dots, a_n) > 0$. Then the following equality holds in $\qq[w](\plgt)$:
$$\sum_{j = 0}^n (-1)^j B_{n, a_0, \dots, \omitvar{a_j}, \dots, a_n}(v) = 0$$
\item Suppose $a_0, \dots, a_n$ are in general position with $k^{-}(a_0, \dots, a_n) = 0$. Then in $\qq[w](\plgt)$:
$$ \sum_{j = 0}^n (-1)^j B_{n, a_0, \dots, \omitvar{a_j}, \dots, a_n}(v) = \begin{cases} \signdet(a_1, \dots, a_n) &\text{ if } v \in L \\ 0 & \text{ otherwise}\end{cases} $$
\end{itemize}
\end{general}

\begin{proof}
If $\rank(a_0, \dots, a_n) \leq n-1$ then the equality is trivial as all terms in the left-hand side vanish. Suppose now that $\rank(a_0, \dots, a_n) = n$. We combine the results from lemma \ref{lemmabernoullicones} with the results from Theorem \ref{theoremcocyclecones}. Denote as before $\eps_j = (-1)^j \signdet(a_0, \dots, \omitvar{a_j}, \dots, a_n)$. From lemma \ref{lemmabernoullicones} we get the equality between rational functions in $\qq[w](\plgt)$:
$$\sum_{j = 0}^n (-1)^j B_{n, a_0, \dots, \omitvar{a_j}, \dots, a_n}(v) = h_0\left(\sum_{j = 0}^n \eps_j\cdual(a_0, \dots, \omitvar{a_j}, \dots, a_n), v\right)$$
If $k^{-}(a_0, \dots, a_n) > 0$ then $\sum_{j = 0}^n \eps_j\cdual( a_0, \dots, \omitvar{a_j}, \dots, a_n) \in \spanqwedges$ by Theorem \ref{theoremcocyclecones}, therefore the right-hand side vanishes. If $a_0, \dots, a_n$ are in general position with $k^{-}(a_0, \dots, a_n) = 0$ then $\sum_{j = 0}^n \eps_j\cdual( a_0, \dots, \omitvar{a_j}, \dots, a_n) - \eps_0 \dirac \in \spanwedges$, therefore:
$$\sum_{j = 0}^n (-1)^j B_{n, a_0, \dots, \omitvar{a_j}, \dots, a_n}(v) = \eps_0 h_0(\dirac, v) =\begin{cases} \eps_0 &\text{ if } v \in L \\ 0 & \text{ otherwise} \end{cases} $$
which proves the claim.
\end{proof}

Unfortunately, the case where $a_0, \dots, a_n$ span $\Lambda$ yet are not in general position with $k^{-}(a_0, \dots, a_n) = 0$ doesn't give good results. Indeed, already for $n = 2$, if $a, b$ are linearly independent with $\det(a,b) = \pm 1$, then:
$$B_{2, -a, b}(0)(w, \plgt) - B_{2, a, b}(0)(w, \plgt) +B_{2, a, -a}(0)(w, \plgt) = \signdet(a,b)\left(\frac{w}{\plgt(\beta)} - \frac12\right) \neq 0 $$
where $\alpha, \beta$ is the primitive positive dual family to $a, b$ in $L$. The right-hand side depends on $b$ and there is no possible value we could have chosen by convention for $B_{2, a, -a}$ which would have given:
$$ B_{2, -a, b}(v)(w, \plgt) - B_{2, a, b}(v)(w, \plgt) + B_{2, a, -a}(v)(w, \plgt) = 0$$
for any non-zero primitive $b \in \Lambda$. 
 
We end this section by showing that in dimension $2$ we may use Theorem \ref{theoremcocyclecones} in conjonction with the Solomon-Hu pairing to recover formula (5.16) in a recent article by Sharifi and Venkatesh \cite{VS} which is a key ingredient in the lifting of a cocycle carried out in [\hspace{1sp}\cite{VS}, Proposition 5.4.1]. Let us write $\theta_{SV}(l_1, l_2)$ for the function $\theta_L(l_1, l_2)$ they define in \symbolparagraph 5.3.1 to avoid confusion with our notations. When $l_1, l_2 \in S^1$ with $\epsilon_{12} = \det(l_1, l_2) \neq 0$ the function $\theta_{SV}(l_1, l_2)$ may be expressed using [\hspace{1sp}\cite{VS}, Lemma 5.3.1] in our notations as:
\begin{equation}
\theta_{SV}(l_1, l_2) = \frac{1-\epsilon_{12}}{2} + \epsilon_{12} h(\cdual(a_1, a_2), 0)(-x)
\end{equation}
where $a_i = \langle \cdot, l_i \rangle$ and $x$ is the linear form corresponding to the formal coordinates $(u_1, u_2)$ defined in [\hspace{1sp}\cite{VS}, \symbolparagraph 5.3.1]. Let us consider $\kappa(l_1, l_2, l_3) = \theta_{SV}(l_1, l_2) + \theta_{SV}(l_2, l_3) -\theta_{SV}(l_1, l_3)$. We wish to show that if $l_1, l_2, l_3$ are in general position then $\kappa(l_1, l_2, l_3) = \delta_{SV}(l_1, l_2, l_3)$ where the function $\delta_{SV}(l_1, l_2, l_3)$ is defined to be $0$ if $l_2$ lies on the counterclockwise portion of $S^1$ joining $l_1$ to $l_3$ and $1$ otherwise. To achieve this we must translate our notion of configurations based on the value of $k^{-}(a_1, a_2, a_3)$ in terms of the function $\delta_{SV}$. Let us then rewrite $\kappa(l_1, l_2, l_3)$ when $l_1, l_2, l_3$ are in general position as: 
$$ \kappa(l_1, l_2, l_3) = \frac{1-\epsilon_{12}-\epsilon_{23}+\epsilon_{13}}{2} + h(\epsilon_{12}\cdual(a_1, a_2) + \epsilon_{23}\cdual(a_2, a_3) - \epsilon_{13}\cdual(a_1, a_3), 0)(-x).$$
where $\epsilon_{ij} = \det(l_i, l_j) = \det(a_i, a_j)$. It follows from Theorem \ref{theoremcocyclecones} that:
$$ \epsilon_{12}\cdual(a_1, a_2) + \epsilon_{23}\cdual(a_2, a_3) - \epsilon_{13}\cdual(a_1, a_3) \equiv \epsilon_{12} q(k^{-}(a_1, a_2, a_3)) \dirac \mod \spanwedges$$
where $q(k) = 0$ if $k > 0$ and $q(0) = 1$ and $\dirac$ is the Dirac function at $0$. Using the linearity of $h$ and its vanishing on $\spanwedges$ we get:
\begin{equation}\label{eqkappa}
\kappa(l_1, l_2, l_3) = \frac{1-\epsilon_{12}-\epsilon_{23}+\epsilon_{13}}{2} + \epsilon_{12} q(k^{-}(a_1, a_2, a_3))
\end{equation}
since $h(\delta, 0)(-x) = 1$. To compute the last term in \refp{eqkappa} we consider $\sum_{j = 1}^3 \lambda_j a_j = 0$ the standard non-trivial relation among $a_1, a_2, a_3$ (see Definition \ref{definitionrelation}) for which it is also true that $\sum_{j = 1}^3 \lambda_j l_j = 0$. Put $\kappa_1(l_1, l_2, l_3) = (1-\epsilon_{12} - \epsilon_{23} + \epsilon_{13})/2$ and $\kappa_2(l_1, l_2, l_3) = \epsilon_{12} q(k^{-}(a_1, a_2, a_3))$ so that $\kappa = \kappa_1 + \kappa_2$. Then we treat $8$ cases separately according to the signs of the $\epsilon_{ij}$'s and gather the results in the following table:

\begin{center}
\begin{tabular}{| C | C | C | C | C | C | C | C | C | C | C |}
\hline
\epsilon_{12} & \epsilon_{13} & \epsilon_{23} & \lambda_1 & \lambda_2 & \lambda_3 & k^{-} & \kappa_1 & \kappa_2 & \kappa & \delta_{SV} \\
\hline
+ & + & + & + & - & + & 1 & 0 & 0 & 0 & 0 \\
\hline
+ & + & - & + & + & - & 1 & 1 & 0 & 1 & 1\\
\hline
+ & - & + & + & + & + & 0 & -1 & 1 & 0 & 0 \\
\hline
+ & - & - & - & + & + & 1 & 0 & 0 & 0 & 0 \\
\hline
- & + & + & - & + & + & 1 & 1 & 0 & 1 & 1 \\
\hline
- & + & - & + & + & + & 0 & 2 & -1 & 1 & 1 \\
\hline
- & - & + & + & + & - & 1 & 0 & 0 & 0 & 0 \\
\hline
- & - & - & + & - & + & 1 & 1 & 0 & 1 & 1\\
\hline
\end{tabular}
\end{center}

This completes the proof of [\hspace{1sp}\cite{VS}, formula (5.16)] when $l_1, l_2, l_3$ are in general position. A similar formula holds for the coefficient of degree $0$ in the expansion of $\theta_{SV}(l_1, l_2)$ near the origin and it may be obtained using Proposition \ref{propositionbernoullicocycle} for $n = 2$. Interesting future work in this direction would be to use Theorem \ref{theoremcocyclecones} in dimension $n \geq 3$ to prove analogues of [\hspace{1sp}\cite{VS}, formula (5.16)] for higher degree analogues of the function $\theta_{SV}(l_1, l_2)$ which appear in recent work by Xu (see \cite{Xu}).

\subsection{The cocycle $\phi_{n,a}$ for unit groups}\label{sectionapplication}

We now define the cocycle properly on specific subgroups of $\slnz{n}$. Fix $a \in \Lambda = \homlz$ a primitive linear form on $L$. We consider subgroups $U$ of $\slnz{n}$ satisfying the following property: $\forall\,m \geq 2, \forall\,g_1, \dots, g_m \in U, \forall \mu_1, \dots, \mu_m \in \zz_{\geq 0}$,
\begin{equation}\label{conditioncocycle}
\sum_{j = 1}^m \mu_j g_j \cdot a = 0 \Rightarrow \mu_1 = 0, \dots, \mu_m = 0
\end{equation}
This property guarantees that we avoid families of vectors in bad position (see \refp{badpositioncondition}) in what follows and we give examples of such groups down below. Recall the definition of the two functions:
$$\psi_{n,a}  := \begin{cases}\slnz{n}^{n-2} & \to \mathcal{F}(V/L \times \cc \times \homlc, \cc) \\ (g_1, \dots, g_{n-2}) & \to \left((v, w, \plgt) \to G_{n-2, a, g_1\cdot a, \dots, (g_1\dots g_{n-2}) \cdot a}(v)(w, \plgt)\right)\end{cases}$$
$$\phi_{n,a}  := \begin{cases} \slnz{n}^{n-1} & \to \mathcal{F}(V/L, \qq[w](\plgt)) \\ (g_1, \dots, g_{n-1}) & \to  B_{n, a, g_1\cdot a, (g_1g_2)\cdot a, \dots, (g_1\dots g_{n-1})\cdot a}(v)(w, \plgt)\end{cases}$$
It follows from Proposition \ref{propositionbernoullicocycle} that $\phi_{n, a}$ is an additive $(n-1)$-cocycle for $U$ as:
$$ \cohomd \phi_{n,a}(g_1, \dots, g_{n+1}) = 0$$
for any $g_1, \dots, g_{n+1} \in U$. Furthermore, from Theorem \ref{theoremmodulargeometricgrintroduction} we deduce that the multiplicative cocycle $\exp(2i\pi\phi_{n, a})$ is partially split by $\psi_{n, a}$. Namely, when $a, g_1 \cdot a, \dots, (g_1 \dots g_{n-1})\cdot a$ are linearly independent, the multiplicative coboundary of $\psi_{n, a}$ is given by:
$$ \cohomdx \psi_{n,a}(g_1, \dots, g_n) = \exp(2i\pi \phi_{n,a}(g_1, \dots, g_n))$$
We expect this splitting property to be true more generally, as we expect that the modular property:
$$ \prod_{j = 1}^n G_{n-2, (a_k)_{k \neq j}}(v)(w, \plgt)^{(-1)^{j+1}} = \exp(2i\pi B_{n, a_1, \dots, a_n}(v)(w, \plgt))$$
holds whenever $a_1, \dots, a_n$ are not \refp{badpositioncondition}, but the strategy of proof we have in mind to extend the domain of validity of Theorem \ref{theoremmodulargeometricgrintroduction} would probably carry us too far from the matter at hand. It should be noted that there are two easier configurations of the $a_i$'s for which the result is true for purely cohomological reasons, and this already covers all cases for $n =2, 3$. In a companion paper in this series, we will focus more on this splitting relation and show that it indeed holds for the subgroups $U$ of $\slnz{n}$ satisfying \refp{conditioncocycle}.

We now give several examples of such groups $U$. A simple example is given by the subgroup of $\slnz{n}$ consisting of matrices which stabilise the set $\{ \alpha \in L \setseparator a(\alpha) \geq 0\}$. This corresponds to the case where linear forms are on the same side of some hyperplane through the origin in [\hspace{1sp}\cite{FDuke}, Lemma 3.9]. It is clear that such a group satisfies condition \refp{conditioncocycle} and that it is isomorphic to the subgroup:
$$U^0 := \left\{ g \in \slnz{n} \setseparator g = \begin{pmatrix} 1 & \ast & \dots & \ast \\ 0 & \ast & \dots & \ast \\ \vdots & \vdots & \vdots & \vdots \\ 0 & \ast & \dots & \ast \end{pmatrix} \right\} $$
of $\slnz{n}$. This simple example is not very interesting for us as it reduces $\phi_{n,a}$ to a $(n-1)$-cocycle for $\slnz{n-1}$ and we now give another set of examples which motivated this work. Consider $\kk$ a number field of degree $n \geq 2$ with at least one real embedding $\sigma_{\rr}$. Denote by $\ok$ the ring of integers of $\kk$ and $\units$ the group of units of $\ok$. By Dirichlet's unit theorem, $\units$ is a free abelian group of rank $r_1 + r_2 -1$ where $r_1$ is the number of real embeddings of $\kk$ and $r_2$ is the number of complex places of $\kk$. Consider now $\mathcal{U}$ a subgroup of $\units$ of full rank such that $\forall \eps \in \mathcal{U}, \sigma_{\rr}(\eps) > 0$. In most examples, $\mathcal{U}$ will be the group of totally positive units $\opck$ of $\ok$. Suppose that $L$ is a lattice of rank $n$ in $\kk$ which is stable under multiplication by elements of $\mathcal{U}$. Fix a $\zz$-basis $B = [e_1, \dots, e_n]$ of $L$. Then $\mathcal{U}$ may be identified with a commutative subgroup $U$ of $\slnz{n}$ which satisfies condition \refp{conditioncocycle} for any non-zero primitive linear form $a \in \Lambda$. Indeed, for any integer $m \geq 2$, if $\eps_1, \dots, \eps_m \in \mathcal{U} \simeq U$ and $\mu_1 \geq 0, \dots, \mu_m \geq 0$ then since:
$$ \left(\sum_{j = 1}^m \mu_j (\eps_j \cdot a)\right)(\cdot)= a\left(\sum_{j = m}^n \mu_j \eps_j^{-1} \cdot \right)$$
it follows that
$$ \sum_{j = 1}^m \mu_j (\eps_j \cdot a) = 0 \Rightarrow a\left(\sum_{j = 1}^m \mu_j \eps_j^{-1} \cdot \right) = 0 \Rightarrow \sum_{j = 1}^m \mu_j \eps_j^{-1} = 0 \Rightarrow \sum_{j = 1}^m \mu_j \sigma_{\rr}(\eps_j^{-1}) = 0$$
as non-zero elements of $\kk$ give bijections of $\kk$ by multiplication. The latter expression is a sum of non-negative numbers which is equal to zero, from which we conclude that $\mu_j = 0, \forall\, 1 \leq j \leq m$. Thus $\phi_{n,a}$ is a true $(n-1)$-cocycle for the group $U \simeq \mathcal{U}$. 

The arithmetic applications in the next section make use of the cocycle properties of $\phi_{n,a}$ for these groups of totally positive units to compute partial zeta values at $s = 0$ in totally real number fields following Shintani's method. In the third paper in this series, we give arithmetic application of the collection of $\psi_{n,a}$ functions. Indeed, in \cite{BCG} Bergeron, Charollois and Garc\'ia use the $\psi_{3, a}$ functions to establish a Kronecker limit formula expressing the values of the first derivative of partial zeta functions at $s=0$ in complex cubic fields in terms of the geometric elliptic $\Gamma_{a,b}$ functions. We believe that such results may be extended to higher degree number fields with exactly one complex place using the geometric $G_{n-2}$ functions we introduced in this article and we have successfully produced many examples to support this claim \cite{preprint}.

\section{Application to the computation of partial zeta functions at $s=0$ for totally real number fields}\label{sectionshintani}

In this last section we express partial zeta values in totally real number fields at $s=0$ in terms of the Bernoulli rational functions $B_{n, a_1, \dots, a_n}$. This gives a connection between the arithmetic of totally real number fields and cocycles extracted from the multiple elliptic Gamma functions. Our approach will use the tools developed by Shintani in \cite{Shintani}.

\subsection{Shintani's method}

Let $\ff$ be a totally real number field of degree $n$ and denote by $\of$ the ring of integers of $\ff$. Fix $\goth{f}$ an integral ideal of $\of$ which will be the finite part of the class field modulus. For any integral ideal $\goth{b}$ coprime to $\goth{f}$ which represents a class in the narrow ray class group mod $\goth{f}$, the partial zeta function attached to $\goth{f}$ and $\goth{b}$ is defined by:
$$\zeta_{\goth{f}}(\goth{b}, s) := \sum_{\goth{a} \sim \goth{b}} \norm{\goth{a}}^{-s}$$
Following Siegel \cite{Siegel} we may rewrite this as:
$$ \zeta_{\goth{f}}(\goth{b}, s) = \norm{\goth{b}}^{-s}\sum_{\mu \in (1+ \goth{f} \goth{b}^{-1})^{+}/\opcf} \norm{\mu}^{-s}$$
where $\opcf$ is the group of totally positive units of $\of$ congruent to $1$ mod $\goth{f}$ and $(1+\goth{f}\goth{b}^{-1})^{+}$ is the set of totally positive elements $v \in \ff$ such that $v-1 \in \goth{f}\goth{b}^{-1}$. The fractional ideal $L = \goth{f} \goth{b}^{-1}$ is a lattice of rank $n$ inside $\ff$. Shintani's strategy to express the values of these partial zeta functions at integers $k \leq 0$ revolves around finding a fundamental domain for the action of $\opcf$ on $(1+L)^{+}$ which can be decomposed in rational polyhedral cones. Shintani then associates to each of these cones a partial zeta function which admits a meromorphic continuation over $\cc - \{1\}$ with values at integers $k \leq 0$ given by specific Bernoulli polynomials. The zeta function associated to a rational polyhedral cone $C \subset \ff^{+} \cup \{0\}$ where $\ff^{+}$ is the set of totally positive elements in $\ff$ and to a vector $v \in \ff/L - \{0\}$:
$$\zeta(C, L, v, s) := \sum_{\mu \in C \cap (v+L)} \norm{\mu}^{-s}$$
We now prove that the values at $s= 0$ of these zeta functions associated to rational polyhedral cones may be expressed as linear combination of the Bernoulli rational functions $B_{n, a_1, \dots, a_n}$.

\begin{general}{Proposition}\label{propositionshintaniatzero}
Let $\alpha_1, \dots, \alpha_m$ be $m$ linearly independent primitive vectors in $L$ where $1 \leq m \leq n$. Suppose that the cone $C = \sum_{j = 1}^m \qq_{\geq 0} \alpha_j$ is included  in $\ff^{+} \cup \{0\}$. Fix $v \in \ff/L - \{0\}$ Denote by $\sigma_1, \dots, \sigma_n$ the real embeddings of $\ff$. Then:
$$\zeta(C, L, v, 0) = \frac1n \sum_{k = 1}^n h_0(C, v)(0, -\sigma_k) $$  
In particular, if $n = m$ and $\cdual(a_1, \dots, a_n) = C$ with $\signdet(a_1, \dots, a_n) = 1$ then:
$$\zeta(C, L, v, 0) = \frac1n \sum_{k = 1}^n B_{n, a_1, \dots, a_n}(v)(0, -\sigma_k) $$  
Furthermore, this number is the trace of an element in $\ff$ and therefore belongs to $\qq$.
\end{general}

\begin{proof}
We follow closely Shintani's original proof of [\hspace{1sp}\cite{Shintani}, Theorem 1]. The norm on $v+L$ is a product of affine linear forms with positive coefficients. Namely, we fix a $\zz$-basis $B = [e_1, \dots, e_n]$ of $L$ consisting of positive elements, and we fix an ordering on the real embeddings $\sigma_1, \dots, \sigma_n$ of $\ff$. If we write $l_{jk} = \sigma_{j}(e_k) > 0$ then for any $\mu = \sum_{k = 1}^n \mu_k e_k \in v + L$ we get:
$$ \norm{\mu} = \prod_{j = 1}^n \sigma_j(\mu) = \prod_{j = 1}^n \sum_{k = 1}^n l_{jk} \mu_k$$
For the rest of this section, $\Gamma(s)$ will denote Euler's $\Gamma$ function. Following Shintani's original proof \cite{Shintani} we get for any rational polyhedral cone $C \subset \ff^{+}$:
$$\Gamma(s)^n\zeta(C, L, v, s) = \int_0^{+\infty} \dots \int_0^{+\infty}h(C, v)\left(-\sum_{j = 1}^nu_j\sigma_j\right)(u_1\dots u_n)^{s-1}du_1\dots du_n  $$
Consider Shintani's decomposition of the positive orthant $\rr_{>0}^n$ given by the sets $D_k = \{ (u_1, \dots, u_n) \in \rr_{\geq 0}^n \setseparator u_i \leq u_k, \forall\, 1 \leq i \leq n\}$. We write the above integrals as integrals over the sets $D_k$ which we handle separately:
$$\Gamma(s)^n\zeta(C, L, v, s) = \sum_{k = 1}^n\int_{D_k}h(C, v)\left(-\sum_{j = 1}^nu_j\sigma_j\right)(u_1\dots u_n)^{s-1}du_1\dots du_n  $$
Let us now introduce variables $t, x_1, \dots, x_n$ such that $t = u_k \in [0, \infty)$ and $tx_i = u_i$ for $i \neq k$ with $x_i \in [0, 1]$. For convenience write $d_k(x, s) = \prod_{j \neq k} x_j^{s-1} dx_j$. Then by a change of variables we get:
$$ \Gamma(s)^n\zeta(C, L, v, s) = \sum_{k = 1}^n \int_{0}^{+\infty} t^{ns-1} \int_0^1 \dots \int_0^1 h(C, v)\left(-t\left(\sigma_k + \sum_{\substack{j = 1 \\ j \neq k}}^nx_j\sigma_j\right)\right)dt \dkxs $$
Let us  now isolate the integral over $t$. Define for $y \in \myhom{\zz}{L}{\rr}$ such that $y(\ff^+) \subset \rr_{>0}$:
$$\chi(s, y) := \frac{1}{\Gamma(s)} \int_0^{+\infty} h(C, v)(-t.y)t^{ns-1}dt$$
It follows from lemma 3.1 in \cite{Colmez} that at $s=0$:
$$\chi(0, y) = \frac{1}{n} h_0(C, v)(0, -y) $$
(see Definition \ref{definitionhzero} for the definition of $h_0$). The value of $\zeta(C, L, v, s)$ at $s = 0$ is then given by:
$$\zeta(C, L, v, 0) = \lim_{s \to 0} \frac{\Gamma(s)^{-(n-1)}}{n} \sum_{k = 1}^n \int_0^1 \dots, \int_0^1 h_0(C, v)\left(0, -\sigma_k - \sum_{j \neq k} x_j \sigma_j\right) \dkxs $$
Now, by definition of $h_0$ we have :
$$h_0(C, v)(0, -\sigma_k - \sum_{j \neq k} x_j \sigma_j) = \frac{\mathrm{polynomial}}{\prod_{l = 1}^m (-\sigma_k - \sum_{j \neq k} x_j \sigma_j)(\alpha_l)} $$
and the denominator does not vanish on the integration domain as each $\alpha_l$ belongs to $\ff^{+}$. Therefore, using a variant of [\hspace{1sp}\cite{Colmez}, Lemma 3.2] we get:
\begin{equation}\label{zetaconevalue}
\zeta(C, L, v, 0) = \frac1n \sum_{k = 1}^n h_0(C, v)(0, -\sigma_k)
\end{equation}
By linearity, this result also applies to open rational polyhedral cones. In the particular case where $n = m$ and $\cdual(a_1, \dots, a_n) = C$ with $\signdet(a_1, \dots, a_n) = 1$ we may rewrite this result using lemma \ref{lemmabernoullicones} as:
\begin{equation}\label{zetaconedualvalue}
\zeta(\cdual(a_1, \dots, a_n), L, v, 0) = \frac1n \sum_{k = 1}^n B_{n, a_1, \dots, a_n}(v)(0, -\sigma_k) 
\end{equation}
It is clear that right-hand side of either \refp{zetaconevalue} or \refp{zetaconedualvalue} is the trace of an algebraic number inside $\ff$ and therefore lies in $\qq$, which was already obtained by Shintani.
\end{proof}

From this point forward, we may take two slightly different approaches to the computation of the full partial zeta functions at $s= 0$, one of them following closely Shintani's original strategy \cite{Shintani} and the other following the ``signed fundamental domain'' strategy in \cite{Diazydiaz}.

\subsection{Shintani's fundamental domain}

The first approach to the computation of partial zeta values at $s=0$ using Proposition \ref{propositionshintaniatzero} follows closely Shintani's original strategy carried out in \cite{Shintani}. We introduce the set:
$$\overline{D} = \{ x \in \ff^{+} \setseparator \forall\,u \in \opcf, \trace((u-1)x) \geq 0 \}$$
which may be written in our notation as $\overline{D} = \cdual(a_u, u \in \opcf)$ where $a_u = \trace((u-1)\cdot)$. In essence, Shintani proved that there is a finite set $E \subset \opcf$ such that in our notation $\overline{D} = \cdual(a_u, u \in E)\cdual(a_{u^{-1}}, u \in E)$ and such that the set $D$ whose indicator function is $\cdual(a_u, u \in E) \prod_{u \in E} (1 - \cdual(-a_{u^{-1}}))$ constitutes a fundamental domain for the action of $\opcf$ on $\ff^{+}$. This already gives the equality:
$$ \zeta_{\goth{f}}(\goth{b}, 0) = \zeta(D, \goth{f}\goth{b}^{-1}, 1_{\ff}, 0)$$
Then the set $D$ may be decomposed as a finite disjoint union of open rational polyhedral cones $C_j, 1 \leq j \leq m$ which all admit a set of linearly independent generators $\{\alpha_{j, 1}, \dots, \alpha_{j, n_j}\} \subset \ff^{+}$ with $1 \leq n_j \leq n$. Putting this together with the results from Proposition \ref{propositionshintaniatzero} gives:
\begin{equation}\label{smallshintanieq}
\zeta_{\goth{f}}(\goth{b}, 0) = \frac1n\sum_{k = 1}^n \sum_{j = 1}^m h_0(C_j, 1_{\ff})(0, -\sigma_k).
\end{equation}
Most of the cones $C_j$ are not full-dimensional (case $n_j < n$) and we would like to obtain a result using only full-dimensional cones. To achieve this, we may use any algebraic manipulations on the right-hand side of \refp{smallshintanieq} using the properties of $h_0$ on $\spancones$ where here $V = L \otimes \qq = \ff$. Indeed, if $f \in \spancones$ is any function congruent to the indicator function of $D$ modulo $\spanwedges$, then we may use the linearity of $h_0$ and its vanishing on $\spanwedges$ to conclude that:
$$ \zeta_{\goth{f}}(\goth{b}, 0) = \frac1n\sum_{k = 1}^n h_0(f, 1_{\ff})(0, -\sigma_k) $$
In particular, we show (see \cite{thesis}) that $\spancones = \spanwedges + \spanconesk{n}$ where $\spanconesk{n}$ the subspace of $\spancones$ spanned by the functions $\cdual(a_1, \dots, a_n)$ for linearly independent $a_1, \dots, a_n$ in $\Lambda$. The proof of such a statement in fact produces an algorithm that expresses any indicator function $\cdual(a_1, \dots, a_m)$ modulo $\spanwedges$ in terms of functions $\cdual(a_{j,1}, \dots, a_{j,n})$ where each of the $a_{j, l}$'s are equal to some $\pm a_k$. This gives the following:

\begin{general}{Proposition}\label{propositionshintani}
There is an integer $N \geq 1$ and there are units $u_{j,l} \in E$ as well as signs $\epsilon_{j, l}, \eta_{j, l} \in \{\pm 1\}$ for $1 \leq j \leq N$ and $1 \leq l \leq n$ such that:
$$\zeta_{\goth{f}}(\goth{b}, 0) = \frac1n \sum_{k = 1}^n \sum_{j = 1}^N \signdet(a_{j, 1}, \dots, a_{j, n})B_{n, a_{j,1}, \dots, a_{j,n}}(1_{\ff})(0, -\sigma_k) $$
where $a_{j, l} = \epsilon_{j,l}\trace((u_{j, l}^{\eta_{j, l}}-1)\cdot)$. Furthermore this decomposition is explicitly computable.
\end{general}

\noindent\textbf{Example:} Consider the real quadratic field $\ff = \qq(\sqrt{19})$. Fix $\goth{f} = (13)$ and $\goth{b} = (1)$. Denote by $\sigma_1, \sigma_2$ the real embeddings of $\ff$. The group of totally positive units congruent to $1 \mod \goth{f}$ is generated by $\eps = 170 + 39\sqrt{19}$. A possible totally positive basis for $L = \goth{f}\goth{b}^{-1}$ is given by $B = [13, 65 + 13\sqrt{19}]$. In this situation we have $D = \cdual(a_1)(1-\cdual(a_{-1}))$ where $a_j = \trace((\eps^j-1)\cdot)$. Reducing modulo $\spanwedges$ gives $D \equiv -\cdual(a_1, a_{-1}) \mod \spanwedges$. Therefore:
$$ \zeta_{\goth{f}}(\goth{b}, 0) = -\frac12 \sum_{k = 1}^2 \signdet(a_1, a_{-1})B_{2, a_{1}, a_{-1}}(1_{\ff})(0, -\sigma_k) $$
We compute explicitly in the basis $B$ the coordinates of the linear forms $a_1$ and $a_{-1}$:
$$a_1 = 338. [13, 122], ~~ a_{-1} = 338. [13, 8] $$
where we have factored in the $\mathrm{gcd}$'s of the coefficients. We carried out the computations using the computer software Pari/GP \cite{parigp} and found:
$$ -\frac12 \sum_{k = 1}^2 \signdet(a_1, a_{-1})B_{2, a_{1}, a_{-1}}(1_{\ff})(0, -\sigma_k) = \trace\left(\frac{33}{104} \right) = \frac{33}{52}$$
which we may check is the value of $\zeta_{(13)}((1), 0)$ using for instance Pari/GP's \textbf{bnrL1} command. \bigskip

Proposition \ref{propositionshintani} is already a great way to express the partial zeta values in terms of these Bernoulli polynomials $B_{n, a_1, \dots, a_n}$ related to the multiple elliptic Gamma functions. However, the explicit computation of this decomposition is quite tedious in general. In addition, the linear forms involved are not quite of the form we hoped for following the discussion at the end of section \ref{sectionapplication}. Indeed, we would like to evaluate the Bernoulli polynomials on a cycle whose shape would resemble the cycle described in [\hspace{1sp}\cite{CDG}, section 2.6] which corresponds to the sign fundamental domain decomposition in \cite{Diazydiaz}. This is the second approach which we carry out in the next section.

\subsection{Signed fundamental domains}

We now prove Theorem \ref{theoremshintani} using the second approach following \cite{Colmez}, \cite{Diazydiaz}. This approach has two benefits: on the one hand the ``signed fundamental domain'' decomposition described in \cite{Diazydiaz} is easier to compute than the Shintani decomposition. On the other hand, the cones involved may be interpreted in terms of algebraic cycles (see \cite{CDG}).

\medskip

\begin{proofbis}{Proof of Theorem \ref{theoremshintani}}
Recall that $\ff$ is a totally real number field of degree $n$, that $\goth{f}$ is an integral ideal in $\of$ and that the integral ideal $\goth{b}$ represents a class in the narrow ray class group at $\goth{f}$. Denote by $\sigma_1, \dots, \sigma_n$ the real embeddings of $\ff$. We wish to prove that:
$$\zeta_{\goth{f}}(\goth{b}, 0) = \frac1n \sum_{k =1}^n \sum_{\rho \in \goth{S}_{n-1}} \nu_{\rho} B_{n, a_{1, \rho}, \dots, a_{n, \rho}}(1_{\ff})(0, -\sigma_k) $$
where the $\nu_{\rho}$'s are signs in $\{-1, 0, +1\}$ and the $a_{j, \rho}$'s are $\qq$-linear forms on $\ff$. To achieve this we will recall the notations from \cite{Diazydiaz} to describe their signed fundamental domain for the action of $\opcf$ on $L = \goth{f} \goth{b}^{-1}$. Fix $\eps_1, \dots, \eps_{n-1}$ be fundamental units for $\opcf$. For any permutation $\rho \in \goth{S}_{n-1}$ and any index $1 \leq i \leq n$ define the elements 
$$f_{i, \rho} = \prod_{j = 1}^{i-1} \eps_{\rho(j)}.$$
We use cones $C_{\rho}$ with generators $f_{1, \rho}, \dots, f_{n, \rho}$, and we now describe how the boundaries of the cones are chosen. Let us consider the set $S$ of permutations $\rho \in \goth{S}_{n-1}$ such that $f_{1, \rho}, \dots, f_{n, \rho}$ are linearly independent over $\qq$. Let us consider the canonical embedding of $\ff$ into $\rr^n$ given by:
$$ \sigma := \begin{cases} \ff & \to \rr^n \\ v & \to (\sigma_1(v), \dots, \sigma_n(v)) \end{cases}$$
For any permutation $\rho \in S$ and any index $1 \leq i \leq n$ we define the signs $\mu_{i, \rho} \in \{-1, +1\}$ by the formula:
$$\mu_{i, \rho} = \frac{\det(\sigma(f_{1, \rho}), \dots, \sigma(f_{i-1, \rho}), e_n, \sigma(f_{i+1, \rho}), \dots, \sigma(f_{n, \rho}))}{\det(\sigma(f_{1, \rho}), \dots, \sigma(f_{n, \rho}))} $$ 
where the determinants are taken in the canonical basis of $\rr^n$ and $e_n = [0, \dots, 0, 1]^{t}$ is the last vector of this basis. Then we may define as in \cite{Diazydiaz} the sets:
$$\rr_{i, \rho} := \begin{cases} [0, +\infty) & \text{ if } \mu_{i, \rho} > 0 \\ (0, +\infty) & \text{ if } \mu_{i, \rho} < 0 \end{cases} $$
Let us then define the cones $C_{\rho} = \sum_{i = 1}^n \rr_{i, \rho} f_{i, \rho}$ for $\rho \in S$ and $C_{\rho} = \sum_{i = 1}^n \rr_{\geq 0} f_{i, \rho}$ if $\rho \not\in S$. We will denote by $c_{\rho}$ the indicator function of the set $C_{\rho}$. It follows from [\hspace{1sp}\cite{Diazydiaz}, Theorem 1] that there are explicit signs $w_{\rho} \in \{-1, 0, 1\}$ for $\rho \in S$ such that the following equality holds in $\spancones$:
$$\sum_{\rho \in S} w_{\rho} \sum_{u \in \opcf} c_{\rho}(u \cdot) = \chi_{\ff^{+}} $$
where $V = \ff \simeq \qq^n$ and $\chi_{\ff^{+}}$ is the indicator function of $\ff^{+}$. This already gives the relation for the partial zeta function:
$$\zeta_{\goth{f}}(\goth{b}, s) = \norm{\goth{b}}^{-s} \sum_{\rho \in S} w_{\rho}\, \zeta(C_{\rho}, L, 1_{\ff}, s) $$
and specialising at $s = 0$ yields by Proposition \ref{propositionshintaniatzero}:
$$\zeta_{\goth{f}}(\goth{b}, 0) = \frac1n \sum_{k = 1}^n \sum_{\rho \in S} w_{\rho}\, h_0(C_{\rho}, 1_{\ff})(0, -\sigma_k).$$
We may rephrase this by describing $C_{\rho}$ in terms of linear forms. Indeed, let us define for any $\rho \in S$ and any index $1 \leq i \leq n$:
$$b_{i, \rho} = \frac{\det(f_{1, \rho}, \dots, \omitvar{f_{i, \rho}}, \dots, f_{n, \rho})}{\det(f_{1, \rho}, \dots, f_{i, \rho}, \dots, f_{n, \rho})} $$
where the determinant is taken relative to any $\qq$-basis of $\ff$. In other words, $b_{i, \rho}$ is the linear form on $\ff$ satisfying $b_{i, \rho}(f_{i, \rho}) = 1$ and $b_{i, \rho}(f_{j, \rho}) = 0$ if $j \neq i$. For any $\rho \in S$ we may split the indices $i \in \{1, \dots, n\}$ into two sets depending on the value of $\mu_{i, \rho}$ by setting $I_{\rho} = \{ 1 \leq i \leq n \setseparator \mu_{i, \rho} > 0\}$ and $J_{\rho} = \{1, \dots, n\} - I_{\rho}$. This gives the following expression for $c_{\rho}$:
$$c_{\rho} = \prod_{i \in I_{\rho}} \cdual(b_{i, \rho}) \prod_{j \in J_{\rho}} (1- \cdual(-b_{j, \rho}))$$
It is then easy to see that we have the following reduction modulo $\spanwedges$:
$$ c_{\rho} \equiv (-1)^{\cardinalshort{J_{\rho}}}\prod_{i \in I_{\rho}} \cdual(b_{i, \rho}) \prod_{j \in J_{\rho}} \cdual(-b_{j, \rho}) \mod \spanwedges.$$
We therefore define $a_{i, \rho}$ to be unique primitive element in $\Lambda = \homlz$ such that $a_{i, \rho} = \lambda_{i, \rho}\mu_{i, \rho}b_{i, \rho}$ with $\lambda_{i, \rho} \in \qq_{>0}$. Using the linearity of $h_0$ and its vanishing on $\spanwedges$ we get the following relation:
$$\zeta_{\goth{f}}(\goth{b}, 0) = \frac1n \sum_{k = 1}^n \sum_{\rho \in S} w_{\rho}(-1)^{\cardinalshort{J_{\rho}}}h_0(\cdual(a_{1,\rho}, \dots, a_{n,\rho}), 1_{\ff})(0, -\sigma_k).$$
Identifying the right-hand side using lemma \ref{lemmabernoullicones} gives the desired result:
$$\zeta_{\goth{f}}(\goth{b}, 0) = \frac1n \sum_{k = 1}^n \sum_{\rho \in \goth{S}_{n-1}} \nu_{\rho}B_{n, a_{1,\rho}, \dots, a_{n,\rho}}(1_{\ff})(0, -\sigma_k) $$
where $\nu_{\rho} = w_{\rho}(-1)^{\cardinalshort{J_{\rho}}}\signdet(a_{1, \rho}, \dots, a_{n, \rho})$ if $\rho \in S$ and $\nu_{\rho} = 0$ otherwise. Note that the signs $w_{\rho}$ and the cones $C_{\rho}$ appearing in the signed fundamental domain decomposition in \cite{Diazydiaz} are explicitly computable and so are the linear forms $a_{i, \rho}$ as well as the sets $J_{\rho}$. Using the explicit definition of the $B_{n, a_{1, \rho}, \dots, a_{n, \rho}}$ we may rewrite the right-hand side as the trace of an element in $\ff$ which implies that $\zeta_{\goth{f}}(\goth{b}, 0) \in \qq$, as was previously known from the theorem of Klingen and Siegel.
\end{proofbis}

This second expression of the partial zeta functions using the signed fundamental domain from \cite{Diazydiaz} is closer to what we had in mind in the discussion carried out in section \ref{sectionapplication} as the cones $C_{\rho}$ admit generators $f_{i, \rho}$ which may be described in terms of an algebraic cycle (see for instance \cite{CDG}). Yet, we would like the linear forms $a_{i, \rho}$ to be described in terms of similar algebraic cycles. Interesting future work would be to somehow construct a ``dual cycle'' to the cycle presented in [\hspace{1sp}\cite{Sczech}, lemma 5] to be evaluated against the Bernoulli polynomials cocycle. Namely, we wish to construct a linear form $a \in \Lambda$ and define the corresponding cones:
$$c'_{\rho} = \cdual(f_{1, \rho} \cdot a, \dots, f_{n, \rho} \cdot a) $$
for any permutation $\rho \in \goth{S}_{n-1}$ where we recall that the action on $\dualsimple{V}$ is defined by $(g \cdot a)(v) = a(g^{-1}v)$. It would hopefully then be possible to find coefficients $\nu'_{\rho} \in \qq$ such that the linear combination $\sum_{\rho \in \goth{S}_{n-1}} \nu'_{\rho}c'_{\rho}$ is congruent to a signed fundamental domain modulo $\spanwedges$. This would then lead to a relation of the form:
$$\zeta_{\goth{f}}(\goth{b}, 0) = \frac1n \sum_{k = 1}^n \sum_{\rho \in \goth{S}_{n-1}} \nu'_{\rho} B_{n, f_{1, \rho} \cdot a, \dots, f_{n, \rho} \cdot a}(1_{\ff})(0, -\sigma_k)$$
which we may rewrite using the cocycle $\phi_{n,a}$ on the subgroup of $\slnz{n}$ corresponding to the torus $\opcf$ as:
$$\zeta_{\goth{f}}(\goth{b}, 0) = \frac1n \sum_{k = 1}^n \sum_{\rho \in \goth{S}_{n-1}} \nu'_{\rho} \phi_{n, a}(\eps_{\rho(1)}, \dots, \eps_{\rho(n)})(1_{\ff})(0, -\sigma_k).$$
In the real quadratic case we may carry out this last approach by setting $a = \det(1_{\ff}, \cdot)/\det(1_{\ff}, \eps)$ where $\eps$ is a generator for $\opcf$. This gives the expression:
$$\zeta_{\goth{f}}(\goth{b}, 0) = -\frac12 \sum_{k = 1}^2 \phi_{2, a}(\eps)(1_{\ff})(0, -\sigma_k)$$
and therefore the partial zeta values at $s=0$ are given by the evaluation of a partial $1$-cocycle for $\slnz{2}$ against a $1$-cycle arising from $\opcf$.

As a last remark on Theorem \ref{theoremshintani}, we stress that it is a reformulation of Shintani's result, borrowing ideas from \cite{Colmez}, \cite{Diazydiaz} and \cite{CDG}, and that there have been many other approaches to the computation of partial zeta values at non-positive integers in totally real fields (see \cite{Bekki} or \cite{CGS} for instance). We argue that the novelty of our work lies in the extraction of arithmetic quantities such as partial zeta values in totally real fields at $s=0$ from the study of the geometric families of $G_{n-2, a_1, \dots, a_{n-1}}$ functions via the associated geometric families of Bernoulli rational functions $B_{n, a_1, \dots, a_n}$. We will now give two examples of computations of partial zeta values at $s=0$ in real cubic fields following the procedure given in the proof of Theorem \ref{theoremshintani}. All the computations were carried out using the computer software Pari/GP \cite{parigp}.

\subsection{Cubic examples}

\subsubsection{First real cubic example}

We now carry out our procedure to recover an example from \cite{CGS}. Consider the real cubic field $\ff = \qq(z)$ where $z$ is a root of the polynomial $x^3-x^2-4x-1$. Fix $\goth{f} = (5)$ and $\goth{b} = (1)$. Fix the basis $B = [5, 5z + 10, 5z^2 - 5z]$ of $L = \goth{f}\goth{b}^{-1}$. A possible choice of fundamental units for $\opcf$ is given by:
$$\eps_1 = 15z^2 + 25z + 6 ~~ ; ~~ \eps_2 = -15z^2 + 20z + 56 $$
Write $\goth{S} = \{ Id, (12)\}$ and compute as in the proof of Theorem  \ref{theoremshintani} the signs $\mu_{i, \rho}$:

\begin{center} \begin{tabular}{| C | C | C |}
\hline 
& &\\
\mu_{1, \mathrm{Id}}= -1  & \mu_{2, \mathrm{Id}} = +1 & \mu_{3, \mathrm{Id}} = -1 \\
& &\\
\hline 
& &\\
\mu_{1, (12)}= +1 & \mu_{2, (12)} = -1 & \mu_{3, (12)} = +1 \\
& &\\
\hline
\end{tabular}
\end{center}
This readily gives $J_{\mathrm{Id}} = \{1, 3\}$ and $J_{(12)} = \{2\}$. The explicit cones described in the proof of Theorem  \ref{theoremshintani} are:
\begin{align*}
C_{\mathrm{Id}} & = \rr_{>0} 1_{\ff} + \rr_{\geq 0} \eps_1 + \rr_{>0} \eps_1\eps_2 = (1-\cdual(-b_{1, \mathrm{Id}}))\cdual(b_{2, \mathrm{Id}})(1-\cdual(-b_{3, \mathrm{Id}})) \\
C_{(12)} & = \rr_{\geq 0} 1_{\ff} + \rr_{> 0} \eps_2 + \rr_{\geq 0} \eps_1\eps_2 = \cdual(b_{1, (12)})(1-\cdual(-b_{2, (12)}))\cdual(b_{3, (12)})\end{align*}
where the linear forms $b_{i, \rho}$ are given on the basis $B$ by:

\begin{center} \begin{tabular}{| C | C | C |}
\hline 
& &\\
b_{1, \mathrm{Id}}= \frac{1}{7}[35, 22, 114]  & b_{2, \mathrm{Id}} = \frac{1}{7}[0, -10, 29] & b_{3, \mathrm{Id}} = \frac{1}{7}[0, 3, -8] \\
& &\\
\hline 
& &\\
b_{1, (12)}= \frac{1}{97}[485, 302, 1588] & b_{2, (12)} = \frac{1}{97}[0, 10, -29] & b_{3, (12)} = \frac{1}{97}[0, 3, 1] \\
& &\\
\hline
\end{tabular}
\end{center}
The corresponding signs $w_{\rho}$ given in \cite{Diazydiaz} are $w_{\mathrm{Id}} = w_{(12)} = 1$. The complete signs $\nu_{\rho}$ are $\nu_{\mathrm{Id}} = -1$ and $\nu_{(12)} = 1$. Following the proof of Theorem \ref{theoremshintani} we set $a_{i, \mathrm{Id}} = 7\mu_{i, \mathrm{Id}}b_{i, \mathrm{Id}}$ and $a_{i, (12)} = 97\mu_{i, (12)}b_{i, (12)}$ for $1 \leq i \leq 3$. We may then compute using formula \refp{definitionrescaledbernoulli}:
\begin{align*}
R_1 &= -\sum_{k = 1}^3 B_{3, a_{1, \mathrm{Id}}, a_{2, \mathrm{Id}}, a_{3, \mathrm{Id}}}(1_{\ff})(0, -\sigma_k) = \trace_{\ff/\qq}\left(\frac{1975z^2-4525z-1424}{120}\right) = \frac{4489}{60}\\
R_2 &= \sum_{k = 1}^3 B_{3, a_{1, (12)}, a_{2, (12)}, a_{3, (12)}}(1_{\ff})(0, -\sigma_k) =  \trace_{\ff/\qq}\left(\frac{-1975z^2+4525z+1448}{120}\right) = \frac{-4453}{60}
\end{align*}
It follows from Theorem \ref{theoremshintani} that 
$$\zeta_{\goth{f}}(\goth{b}, 0) = \frac13 (R_1 + R_2) = \frac15 $$
which recovers the result in \cite{CGS}.

\subsubsection{Second real cubic example}

We now study an example with a modulus $\goth{f}$ which is not of the form $N\of$ for some rational integer $N \geq 1$. Consider the real cubic field $\ff = \qq(z)$ where $z$ is a root of the polynomial $x^3-x^2-6x+3$. Fix $\goth{f} = (1-z)\of$ the unramified prime ideal above $3$ in $\of$ and $\goth{b} = (1)$. Fix the basis $B = [3, x + 5, x^2 + 2]$ of $L = \goth{f}\goth{b}^{-1}$. A possible choice of fundamental units for $\opcf$ is given by:
$$\eps_1 = -4z^2 + z + 28 ~~ ; ~~ \eps_2 = -3z^2 + 3z + 22 $$
The explicit cones described in the proof of Theorem \ref{theoremshintani} are:
\begin{align*}
C_{\mathrm{Id}} & = \rr_{\geq 0} 1_{\ff} + \rr_{> 0} \eps_1 + \rr_{\geq 0} \eps_1\eps_2 = \cdual(b_{1, \mathrm{Id}})(1-\cdual(-b_{2, \mathrm{Id}}))\cdual(-b_{3, \mathrm{Id}}) \\
C_{(12)} & = \rr_{> 0} 1_{\ff} + \rr_{\geq 0} \eps_2 + \rr_{> 0} \eps_1\eps_2 = (1-\cdual(-b_{1, (12)}))\cdual(b_{2, (12)})(1-\cdual(-b_{3, (12)}))\end{align*}
We directly give the signs $\nu_{\mathrm{Id}} = 1$ and $\nu_{(12)} = -1$ and the expression of the $a_{i, \rho}$'s on the basis $B$ as:

\begin{center} \begin{tabular}{| C | C | C |}
\hline
& & \\
a_{1, \mathrm{Id}}= [108, 280, 349]  & a_{2, \mathrm{Id}} = [0, 25, 13] & a_{3, \mathrm{Id}} = [0, 4, 1] \\
& & \\
\hline 
& & \\
a_{1, (12)}= [-432, -395, -1019] & a_{2, (12)} = [0, 25, 13] & a_{3, (12)} = [0, 1, 1] \\
& & \\
\hline
\end{tabular}
\end{center}
Thus we may compute using formula \refp{definitionrescaledbernoulli} :
\begin{align*}
R_1 &= \sum_{k = 1}^3 B_{3, a_{1, \mathrm{Id}}, a_{2, \mathrm{Id}}, a_{3, \mathrm{Id}}}(1_{\ff})(0, -\sigma_k) = \trace_{\ff/\qq}\left(\frac{3z^2-7}{6}\right) = 3\\
R_2 &= -\sum_{k = 1}^3 B_{3, a_{1, (12)}, a_{2, (12)}, a_{3, (12)}}(1_{\ff})(0, -\sigma_k) =  \trace_{\ff/\qq}\left(\frac{-3z^2+11}{6}\right) = -1
\end{align*}
It follows from Theorem \ref{theoremshintani} that 
$$\zeta_{\goth{f}}(\goth{b}, 0) = \frac13 (R_1 + R_2) = \frac23 $$
which can be verified using Pari/GP's \textbf{bnrL1} command for instance.

\bibliographystyle{alpha}
\bibliography{bibliographie}{}

@misc{BCG,
	author = {Bergeron, Nicolas and Charollois, Pierre and Garc\'ia, Luis E.},
	title = {Elliptic units for complex cubic fields},
	eprint={2311.04110},
	archivePrefix={arXiv},
	journal ={},
	year = {2023},
	howpublished = {https://arxiv.org/abs/2311.04110}
}

@article{CDG,
 author = {Charollois, Pierre and Dasgupta, Samit and Greenberg, Matthew},
 title = {Integral {Eisenstein} cocycles on {{\(\mathrm{GL}_n\)}}. {II}: {Shintani}'s method.},
 fjournal = {Commentarii Mathematici Helvetici},
 journal = {Comment. Math. Helv.},
 issn = {0010-2571},
 volume = {90},
 number = {2},
 pages = {435--477},
 year = {2015},
 language = {English},
 doi = {10.4171/CMH/360},
 zbMATH = {6451266},
 Zbl = {1326.11072}
}

@article{Ruijsenaars,
 author = {Ruijsenaars, S. N. M.},
 title = {First order analytic difference equations and integrable quantum systems},
 fjournal = {Journal of Mathematical Physics},
 journal = {J. Math. Phys.},
 issn = {0022-2488},
 volume = {38},
 number = {2},
 pages = {1069--1146},
 year = {1997},
 language = {English},
 doi = {10.1063/1.531809},
 url = {ir.cwi.nl/pub/2164},
 zbMATH = {1014649},
 Zbl = {0877.39002}
}

@article{FV,
 author = {Felder, Giovanni and Varchenko, Alexander},
 title = {The elliptic gamma function and {{\(\text{SL}(3,\mathbb Z)\ltimes\mathbb Z^3\)}}.},
 fjournal = {Advances in Mathematics},
 journal = {Adv. Math.},
 issn = {0001-8708},
 volume = {156},
 number = {1},
 pages = {44--76},
 year = {2000},
 language = {English},
 doi = {10.1006/aima.2000.1951},
 zbMATH = {1579912},
 Zbl = {1038.11029}
}

@article{FDuke,
 author = {Felder, Giovanni and Henriques, Andr{\'e} and Rossi, Carlo A. and Zhu, Chenchang},
 title = {A gerbe for the elliptic gamma function},
 fjournal = {Duke Mathematical Journal},
 journal = {Duke Math. J.},
 issn = {0012-7094},
 volume = {141},
 number = {1},
 pages = {1--74},
 year = {2008},
 language = {English},
 doi = {10.1215/S0012-7094-08-14111-0},
 url = {ora.ox.ac.uk/objects/uuid:8c728bbe-3e8c-4b4d-8c28-d68c5220a735},
 zbMATH = {5224873},
 Zbl = {1130.33010}
}

@article{Nishizawa,
 author = {Nishizawa, Michitomo},
 title = {An elliptic analogue of the multiple gamma function},
 fjournal = {Journal of Physics A: Mathematical and General},
 journal = {J. Phys. A, Math. Gen.},
 issn = {0305-4470},
 volume = {34},
 number = {36},
 pages = {7411--7421},
 year = {2001},
 language = {English},
 doi = {10.1088/0305-4470/34/36/320},
 zbMATH = {1695708},
 Zbl = {0993.33016}
}

@article{Narukawa,
 author = {Narukawa, Atsushi},
 title = {The modular properties and the integral representations of the multiple elliptic gamma functions},
 fjournal = {Advances in Mathematics},
 journal = {Adv. Math.},
 issn = {0001-8708},
 volume = {189},
 number = {2},
 pages = {247--267},
 year = {2004},
 language = {English},
 doi = {10.1016/j.aim.2003.11.009},
 zbMATH = {2136504},
 Zbl = {1077.33024}
}

@article{Sczech,
 author = {Sczech, Robert},
 title = {Eisenstein group cocycles for {{\(\text{GL}_ n\)}} and values of {{\(L\)}}- functions},
 fjournal = {Inventiones Mathematicae},
 journal = {Invent. Math.},
 issn = {0020-9910},
 volume = {113},
 number = {3},
 pages = {581--616},
 year = {1993},
 language = {English},
 doi = {10.1007/BF01244319},
 url = {https://eudml.org/doc/144138},
 zbMATH = {549667},
 Zbl = {0809.11029}
}

@article{Winding,
 author = {Winding, Jacob},
 title = {Multiple elliptic gamma functions associated to cones},
 fjournal = {Advances in Mathematics},
 journal = {Adv. Math.},
 issn = {0001-8708},
 volume = {325},
 pages = {56--86},
 year = {2018},
 language = {English},
 doi = {10.1016/j.aim.2017.11.022},
 zbMATH = {6824535},
 Zbl = {1390.33009}
}

@manual{parigp,
	organization = "{The PARI~Group}",
	title = "{PARI/GP version \texttt{2.15.4}}",
	year = 2023,
	address = "Univ. Bordeaux",
	note = "available from http://pari.math.u-bordeaux.fr/"
}

@article{Colmez,
 author = {Colmez, Pierre},
 title = {Residue at {{\(s=1\)}} of {{\(p\)}}-adic zeta functions},
 fjournal = {Inventiones Mathematicae},
 journal = {Invent. Math.},
 issn = {0020-9910},
 volume = {91},
 number = {2},
 pages = {371--389},
 year = {1988},
 language = {French},
 doi = {10.1007/BF01389373},
 url = {https://eudml.org/doc/143545},
 zbMATH = {4061355},
 Zbl = {0651.12010}
}

@article{Diazydiaz,
 author = {Diaz y Diaz, Francisco and Friedman, Eduardo},
 title = {Signed fundamental domains for totally real number fields},
 fjournal = {Proceedings of the London Mathematical Society. Third Series},
 journal = {Proc. Lond. Math. Soc. (3)},
 issn = {0024-6115},
 volume = {108},
 number = {4},
 pages = {965--988},
 year = {2014},
 language = {English},
 doi = {10.1112/plms/pdt025},
 zbMATH = {6298167},
 Zbl = {1325.11117}
}

@book{BarvinokConvexity,
 author = {Barvinok, Alexander},
 title = {A course in convexity},
 fseries = {Graduate Studies in Mathematics},
 series = {Grad. Stud. Math.},
 issn = {1065-7338},
 volume = {54},
 isbn = {0-8218-2968-8},
 year = {2002},
 publisher = {Providence, RI: American Mathematical Society (AMS)},
 language = {English},
 zbMATH = {1860211},
 Zbl = {1014.52001}
}

@article{SolomonHu,
 author = {Hu, Shubin and Solomon, David},
 title = {Properties of higher-dimensional {Shintani} generating functions and cocycles on {{\(\text{PGL}_ 3(\mathbb{Q})\)}}.},
 fjournal = {Proceedings of the London Mathematical Society. Third Series},
 journal = {Proc. Lond. Math. Soc. (3)},
 issn = {0024-6115},
 volume = {82},
 number = {1},
 pages = {64--88},
 year = {2001},
 language = {English},
 doi = {10.1112/S0024611500012612},
 zbMATH = {1696300},
 Zbl = {1045.11035}
}

@article{RichardHill,
 author = {Hill, Richard M.},
 title = {Shintani cocycles on {{\(\mathrm{GL}_n\)}}},
 fjournal = {Bulletin of the London Mathematical Society},
 journal = {Bull. Lond. Math. Soc.},
 issn = {0024-6093},
 volume = {39},
 number = {6},
 pages = {993--1004},
 year = {2007},
 language = {English},
 doi = {10.1112/blms/bdm099},
 zbMATH = {5250870},
 Zbl = {1192.11030}
}

@misc{Barnes,
 author = {Barnes, E. W.},
 title = {On the theory of the multiple {Gamma} function.},
 year = {1904},
 language = {English},
 howpublished = {Trans. Cambridge Philos. Soc. 19, 374-425},
 zbMATH = {2653713},
 JFM = {35.0462.01}
}

@misc{preprint,
 author = {Morain, Pierre L. L.},
 title = {Elliptic units above fields with exactly one complex place},
 year = {2024},
 howpublished = {Preprint, {arXiv}:2406.06094 [math.{NT}]},
 url = {https://arxiv.org/abs/2406.06094},
 arXiv = {arXiv:2406.06094}
}

@misc{thesis,
	author = {Morain, Pierre L. L.},
	title = {Geometric families of multiple elliptic {Gamma} functions and arithmetic applications [{PhD} thesis in progress, {Sorbonne Université}]},
	year = {expected 2026}
}

@misc{secondpaper,
	author = {Morain, Pierre L. L.},
	title = {Geometric families of multiple elliptic {Gamma} functions and arithmetic applications, {II} [{In} progress]},
	year = {2026}
}

@misc{thirdpaper,
	author = {Morain, Pierre L. L.},
	title = {Geometric families of multiple elliptic {Gamma} functions and arithmetic applications, {III} [{In} progress]},
	year = {2026}
}

@book{Siegel,
 author = {Siegel, Carl Ludwig},
 title = {Advanced analytic number theory},
 fseries = {Tata Institute of Fundamental Research. Studies in Mathematics},
 series = {Tata Inst. Fundam. Res., Stud. Math.},
 volume = {9},
 year = {1980},
 publisher = {Bombay: Tata Institute of Fundamental Research},
 language = {English},
 zbMATH = {3751041},
 Zbl = {0478.10001}
}

@article{Shintani,
 author = {Shintani, Takuro},
 title = {On evaluation of zeta functions of totally real algebraic number fields at non-positive integers},
 fjournal = {Journal of the Faculty of Science. Section I A},
 journal = {J. Fac. Sci., Univ. Tokyo, Sect. I A},
 issn = {0040-8980},
 volume = {23},
 pages = {393--417},
 year = {1976},
 language = {English},
 zbMATH = {3544161},
 Zbl = {0349.12007}
}

@article{DarmonPozziVonk,
 author = {Darmon, Henri and Pozzi, Alice and Vonk, Jan},
 title = {The values of the {Dedekind}-{Rademacher} cocycle at real multiplication points},
 fjournal = {Journal of the European Mathematical Society (JEMS)},
 journal = {J. Eur. Math. Soc. (JEMS)},
 issn = {1435-9855},
 volume = {26},
 number = {10},
 pages = {3987--4032},
 year = {2024},
 language = {English},
 doi = {10.4171/JEMS/1344},
 zbMATH = {7891436}
}

@incollection{CGS,
 author = {Chinta, Gautam and Gunnells, Paul E. and Sczech, Robert},
 title = {Computing special values of partial zeta-functions},
 booktitle = {Algorithmic number theory. 4th international symposium. ANTS-IV, Leiden, the Netherlands, July 2--7, 2000. Proceedings},
 isbn = {3-540-67695-3},
 pages = {247--256},
 year = {2000},
 publisher = {Berlin: Springer},
 language = {English},
 doi = {10.1007/10722028_13},
 zbMATH = {1643929},
 Zbl = {1032.11049}
}

@article{Spiridonov,
 author = {Spiridonov, V. P.},
 title = {Theta hypergeometric integrals},
 fjournal = {St. Petersburg Mathematical Journal},
 journal = {St. Petersbg. Math. J.},
 issn = {1061-0022},
 volume = {15},
 number = {6},
 pages = {929--967},
 year = {2004},
 language = {English},
 doi = {10.1090/S1061-0022-04-00839-8},
 zbMATH = {2165587},
 Zbl = {1071.33011}
}

@article{VS,
 author = {Sharifi, Romyar and Venkatesh, Akshay},
 title = {Eisenstein cocycles in motivic cohomology},
 fjournal = {Compositio Mathematica},
 journal = {Compos. Math.},
 issn = {0010-437X},
 volume = {160},
 number = {10},
 pages = {2407--2479},
 year = {2024},
 language = {English},
 doi = {10.1112/S0010437X24007322},
 zbMATH = {7940509}
}

@misc{Xu,
 author = {Xu, Peter},
 title = {Symbols for toric {Eisenstein} cocycles and arithmetic applications},
 year = {2025},
 howpublished = {Preprint, {arXiv}:2402.00294 [math.{NT}]},
 url = {https://arxiv.org/abs/2402.00294},
 arXiv = {arXiv:2402.00294}
}

@article{Bekki,
 author = {Bekki, Hohto},
 title = {On the conical zeta values and the {Dedekind} zeta values for totally real fields},
 fjournal = {Acta Arithmetica},
 journal = {Acta Arith.},
 issn = {0065-1036},
 volume = {216},
 number = {2},
 pages = {177--196},
 year = {2024},
 language = {English},
 doi = {10.4064/aa231026-21-5},
 zbMATH = {7958562}
}

\end{document}